\documentclass{gtart}

\input gtmonout
\volumenumber{1}
\volumeyear{1998}
\volumename{The Epstein birthday schrift}
\pagenumbers{511}{549}
\received{15 November 1997}\revised{27 November 1998}
\published{30 November 1998}
\papernumber{25}

\usepackage{graphicx}
\usepackage{floatflt}
\usepackage{amssymb,amsmath,amsthm}
\usepackage{psfrag}
\theoremstyle{plain}
\newtheorem{theorem}{Theorem}[section]
\newtheorem{proposition}[theorem]{Proposition}

\theoremstyle{remark}
\newtheorem*{remark}{Remark}

\newcommand{\gozto}{\mapsto}
\newcommand{\newoperator}[2]{\DeclareMathOperator{#1}{#2}}

\newoperator{\AH}{AH}
\newoperator{\AJ}{AJ}
\newoperator{\BS}{BS}               
\newoperator{\GF}{GF}               
\newoperator{\GH}{GH}
\newoperator{\GJ}{GJ}
\newoperator{\Isom}{Isom}
\newoperator{\PSL}{PSL}
\newoperator{\SL}{SL}
\newoperator{\KS}{KS}               
\newoperator{\ML}{ML}
\newoperator{\Nbhd}{Nbhd}
\newoperator{\PL}{PL}
\newoperator{\QF}{QF}
\newoperator{\QJ}{QJ}
\newoperator{\QpH}{Q^\prime H}
\newoperator{\Rad}{Rad}
\newoperator{\Tr}{Tr}
\newoperator{\Vol}{Vol}
\newoperator{\area}{area}
\newoperator{\base}{base}
\newoperator{\codim}{codim}
\newoperator{\const}{const}
\newoperator{\glue}{glue}
\newoperator{\graft}{graft}
\newoperator{\inj}{inj}             
\newoperator{\interior}{int}
\newoperator{\length}{length}
\newoperator{\llength}{\underline length}
\newoperator{\mass}{mass}
\newoperator{\ls}{ls}
\newoperator{\qf}{qf}
\newoperator{\skin}{skin}
\newoperator{\teich}{\textsl{T}}
\newoperator{\unglue}{unglue}
\newoperator{\wb}{wb}               
\newoperator{\window}{window}
\newoperator{\ws}{ws}               
\newoperator{\bom}{\boundary_0(M)}
\newoperator{\dl}{dl}
\newoperator{\dr}{dr}
\newoperator{\Acyl}{Acyl}
\newoperator{\acyl}{Acyl}
\newoperator{\eisen}{\bf Eis}
\DeclareMathOperator{\isom}{isom}
\newcommand{\df}[1]{\emph{#1}}
\newcommand{\hy}{\mathbb{H}}
\newcommand{\integers}{\mathbb{Z}}

\newcommand{\complexes}{\mathbb{C}}
\newcommand{\CH}{\mathbb{CH}}

\newcommand{\CP}{\mathbb{CP}}
\newcommand{\reals}{\mathbb{R}}
\newcommand{\R}{\reals}
\newcommand{\rationals}{\mathbb{Q}}
\newcommand{\euclidean}{\mathbb{E}}
\newcommand{\fund}[1]{\pi_1(#1)}
\newcommand{\set}[1]{\{#1\}}

\newcommand{\Chat}{\hat{\complexes}}

\newcommand{\arrow}{\rightarrow}
\newcommand{\boundary}{\partial}

\newcommand{\cross}{\times}

\newcommand{\degrees}{^\circ}

\newcommand{\union}{\cup}

\newcount{\examplenumber}
\newcommand\nex{\global\advance\examplenumber by 1 \number\examplenumber.}

\begin{document}
\title{Shapes of polyhedra and triangulations of the sphere}
\author{William P Thurston}
\address{Mathematics Department\\University of California at Davis\\
Davis, CA 95616, USA}
\email{wpt@math.ucdavis.edu}

\begin{abstract}
The space of shapes of a polyhedron with 
given total
angles less than $2 \pi$ at each of its $n$ vertices
has a K\"ahler metric, locally isometric to complex hyperbolic space
$\CH^{n-3}$.
The metric is not complete: collisions between vertices take place
a finite distance from a nonsingular point.  The metric completion
is a complex hyperbolic cone-manifold.  In some interesting special cases,
the metric completion is an orbifold.
The concrete description of these spaces of shapes gives information about
the combinatorial classification of
triangulations of the sphere with no more than $6$ triangles at a vertex.
\end{abstract}

\asciiabstract{The space of shapes of a polyhedron with 
given total
angles less than 2\pi at each of its n vertices
has a Kaehler metric, locally isometric to complex hyperbolic space
CH^{n-3}.
The metric is not complete: collisions between vertices take place
a finite distance from a nonsingular point.  The metric completion
is a complex hyperbolic cone-manifold.  In some interesting special cases,
the metric completion is an orbifold.
The concrete description of these spaces of shapes gives information about
the combinatorial classification of
triangulations of the sphere with no more than 6 triangles at a vertex.}

\primaryclass{51M20}\secondaryclass{51F15, 20H15, 57M50}
\keywords{Polyhedra, triangulations, configuration spaces, braid groups,
complex hyperbolic orbifolds}

\maketitle

\begin{figure}[htbp]
\begin{minipage}{.45\textwidth}
\centering
\includegraphics[width=.9\textwidth]{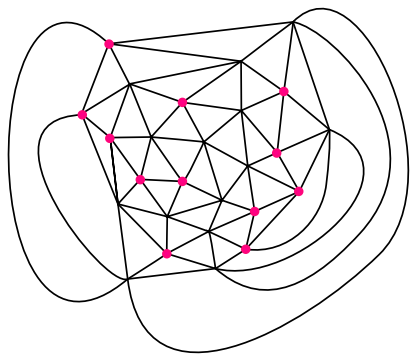}
\vglue -4mm
\end{minipage}
\begin{minipage}{.03\textwidth}\hfil
\end{minipage}
\begin{minipage}{.43\textwidth}
\centering
\caption[A triangulation with 5 or 6 triangles per vertex.]
{The twelve marked vertices of this triangulation of $S^2$ have 
five triangles while all other vertices
have six. Theorem \ref{Polyhedra are lattice points} implies that
the possible triangulations satisfying this condition are parametrized,
up to isomorphism, by 20--tuples of integers up to the action of a group
of integer linear transformations.
}
\label{figure: order 5 and 6 triangulations}
\end{minipage}
\end{figure}

\section*{Introduction}

There are only three completely symmetric triangulations of the sphere:
the tetrahedron, the octahedron and the icosahedron.
However, finer triangulations with good geometric properties 
are often encountered or desired for mathematical,
scientific or technological reasons, for example, the kinds of
triangulations popularized in modern times by
Buckminster Fuller and used for geodesic domes and chemical `Buckyballs'.

There are procedures 
to refine and modify any triangulation of a surface until 
every vertex has either 5, 6 or 7 triangles around it, or with more effort,
so that there are
only 5 or 6 triangles if the surface has positive Euler characteristic,
only 6 triangles if the surface has zero Euler characteristic, or only
6 or 7 triangles if the surface has negative Euler characteristic.
These conditions on triangulations are combinatorial
analogues of metrics of positive, zero or negative curvature.
How systematically can they be understood?

In this paper, we will develop a global theory to
describe all triangulations of the $S^2$
such that each vertex has $6$ or fewer triangles at any vertex.  
Here is one description:

\begin{theorem}[Polyhedra are lattice points]\label{Polyhedra are lattice points}
There is a lattice $L$ in complex Lorentz space $C^{(1,9)}$ and a group
$\Gamma$ of automorphisms, such that triangulations of non-negative
combinatorial curvature are elements of $L_+/\Gamma$, where $L_+$ is
the set of lattice points of positive square-norm.
The projective action of $\Gamma$
on complex hyperbolic space $\CH^9$
(the unit ball in $\complexes^9 \subset \CP^9$)
has quotient of
finite volume.  The square of the norm of a lattice point is the number of
triangles in the triangulation.  \end{theorem}

A triangulation is \df{non-negatively curved} if there are
never more than six triangles at a vertex.
The theorem can be interpreted as describing certain concrete
cut-and-glue constructions for creating triangulations of non-negative
curvature, starting from simple and easily-classified examples.
The constructions are parametrized by choices
of integers, subject to certain geometric constraints.
The fact that $\Gamma$ is a discrete group means
that it is possible to dispense with most of the 
constraints, except for an algebraic condition that a certain
quadratic form is positive: any choice
of integer parameters can be transformed by $\Gamma$ to satisfy the
geometric conditions, and the resulting triangulation is unique. Thus,
the collection of all triangulations can be described either as a quotient
space, in which identifications of the parameters are made algebraically,
or as a fundamental domain (see section \ref{section: construction and fundamental domain}).

\begin{floatingfigure}[.9]{1in}
\psfrag{q}{$\scriptstyle \theta$}
\centering
\includegraphics[width=.7in]{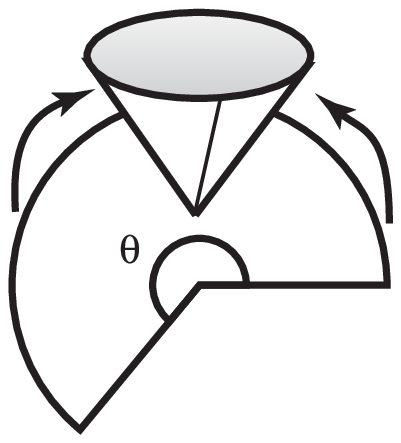}
\small Cone angle $\theta$
\end{floatingfigure}

We will study combinatorial types of triangulations by using a
metric where each triangle is a Euclidean equilateral triangle with
sides of unit length. This metric
is locally Euclidean everywhere except near
vertices that have fewer than $6$ triangles.

It is helpful to consider
these metrics as a special case of metrics on the sphere which are
locally Euclidean except at a finite number of points, 
which have neighborhoods locally modelled on cones. A \df{cone
of cone-angle $\theta$} is a metric space that can be 
formed, if $\theta \le 2 \pi$, from a sector of the Euclidean 
plane between two rays that make an angle $\theta$, by 
gluing the two rays together.  More generally, a cone of 
angle $\theta$ can be formed by taking the universal cover 
of the plane minus $0$, reinserting $0$,  and then 
identifying modulo a transformation that 
``rotates'' by angle $\theta$.  The \df{apex curvature} of a cone 
of cone-angle $\theta$ is $2 \pi - \theta$.

A Euclidean cone metric on a surface satisfies the Gauss--Bonnet theorem,
that is, the sum of the apex curvatures is $2 \pi$ times the Euler characteristic.
This fact can be derived from basic Euclidean geometry
by subdividing the surface into triangles and looking at the sum of angles
of all triangles grouped in two different ways, by triangle or by vertex.
It can also be derived from the usual smooth Gauss--Bonnet formula by
rounding off the points, replacing a tiny neighborhood of each cone point
by a smooth surface (for example part of a small sphere).

\begin{theorem}[Cone metrics form cone manifold]
\label{cone metrics form cone-manifold}
Let $k_1, k_2, \dots, k_n \quad [n > 3]$
be a collection of real numbers in the interval
$(0, 2 \pi)$ whose sum is $4 \pi$.  Then the set of Euclidean cone metrics
on the sphere with cone points of curvature $k_i$ and
of total area $1$ forms a complex hyperbolic manifold, whose
metric completion is a complex hyperbolic cone manifold of finite
volume.
This cone manifold is an orbifold (that is, the quotient space of
a discrete group) if and only if
for any pair $k_i, k_j$ whose sum $s = k_i + k_j$ that
is less than $2 \pi$, either
\begin{description}
\item[(i)] $(2 \pi - s)$ divides $2 \pi$, or
\item[(ii)] $k_i = k_j$, and $(2 \pi - s)/2 = (\pi - k_i)$ divides $2 \pi$.
\end{description}
\end{theorem}

The definition of ``cone-manifold'' in dimensions bigger than
$2$ will be given later.

This turns out to be closely related to work of Picard
(\cite{Picard:hypergeometriques-2}, \cite{Picard:hypergeometriques-1})
and Mostow and Deligne (\cite {Deligne-Mostow:Monodromy},
\cite{Deligne-Mostow:Commensurabilities}, \cite {Deligne-Mostow:Picard}).
Picard discovered
many of the orbifolds; his student LeVavasseur enumerated the class
of groups Picard discovered, and
they were further analyzed by Deligne and Mostow. Mostow discovered that
condition (i) is not always required to obtain an orbifold and that 
(ii) is sufficient when $k_i = k_j$.
However, the geometric interpretations were not apparent in these papers.
It is possible to understand the quotient cone-manifolds quite concretely
in terms of shapes of polyhedra.

A version of this paper has circulated for a number of years as a preprint,
which for a time was circulated as a Geometry Center preprint, and later
revised as part of the xxx mathematics archive.  In view of this
history, some time warp is inevitable: for some parts of this
paper, others may have have done further work that is not here taken into
account.  I would like to thank Derek Holt, Igor Rivin,
Chih-Han Sah and Rich Schwarz for mathematical comments and corrections
that I hope I have taken into account.

\section{Triangulations of a hexagon}
\label{triangulations}

Let $E$ be the standard equilateral triangulation
of $\complexes$ by triangles of unit side length, where $0$ and $1$
are both vertices.
The set $\eisen$ of vertices of $E$ are complex
numbers of the form $m + p \omega $,
where $\omega = 1/2 + \sqrt {-3} / 2$ is a primitive $6$th root of
unity.
These lattice points form a subring of $\complexes$, called
the Eisenstein integers, the ring of algebraic
integers in the quadratic imaginary field $\rationals(\sqrt{-3})$.

\begin{figure}[hbtp]
\begin{minipage}{.01\textwidth}\hfil
\end{minipage}
\begin{minipage}{.30\textwidth}
\includegraphics[width=\textwidth]{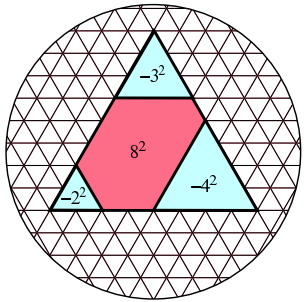}
\end{minipage}
\begin{minipage}{.05\textwidth}\hfil
\end{minipage}
\begin{minipage}{.62\textwidth}
\caption[A hexagon is a big triangle minus three little triangles.]
{An Eisenstein lattice hexagon has the form of a large equilateral
triangle of sidelength $n$,
minus three equilateral triangles that fit inside it of sidelengths
$p_1$, $p_2$ and $p_3$.
An equilateral triangle of sidelength $n$ contains $n^2$ unit equilateral
triangles, so the hexagon has $n^2 - p_1^2 - p_2^2 - p_3^2$ triangles.
}
\end{minipage}
\end{figure}

To warm up, we'll analyze all possible shapes
of Eisenstein lattice hexagons, with vertices in $\eisen$
and sides parallel (in order) with the sides of a standard hexagon.
Note that any such hexagon with $m$ triangles determines a
non-negatively curved
triangulation of the sphere with $2m$ triangles, formed by making
a hexagonal envelope from two copies of the hexagon glued along
the boundary.

If we circumscribe a lattice triangle $T$ about our lattice hexagon $H$, this
gives a description 
\[
H = T \setminus \left ( S_1 \union S_2 \union S_3
\right ) ,
\]
where the $S_i$ are smaller equilateral triangles. If $T$ has sidelength
$n$ and $S_i$ has sidelength $p_i$, then $H$ contains
\begin{equation}
m = n^2 - p_1^2 - p_2^2 - p_3^2
\label{eq: area of hexagon}
\end{equation}
triangles.

All such hexagons are described by four integer parameters,
subject to the 6 inequalities 
\begin{gather*}
p_1 \ge 0 \quad      p_2 \ge 0      \quad p_3 \ge 0 \\
p_1 + p_2 \le n \quad p_2 + p_3 \le n \hss\quad p_3 + p_1 \le n
, \end{gather*}
where strict inequalities give
non-degenerate hexagons; if one or more inequality becomes an equality
then one or more sides of the hexagon shrinks to length 0 and the
`hexagon' becomes a pentagon, quadrilateral or triangle.

\begin{figure}[hbtp]
\centering
\begin{minipage}{.4\textwidth}
\includegraphics[width=\textwidth]{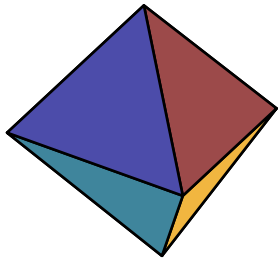}
\end{minipage}
\begin{minipage}{.59\textwidth}
\caption{The space of shapes of hexagons is described by this polyhedron
in hyperbolic 3--space; the faces represent hexagons degenerated to pentagons,
and the edges represent degeneration to quadrilaterals.
  All dihedral angles are $\pi/2$.
The three mid-level vertices are ideal vertices at infinity, and 
represent the three ways that hexagons
can become arbitrarily long and skinny,
while the top and bottom are finite vertices, representing the
two ways that hexagons can degenerate
to equilateral triangles. The polyhedron has hyperbolic
volume $.91596559417\dots$.
\label{figure: hexagon classification}
}
\end{minipage}
\end{figure}
The solutions are elements of the integer lattice inside a convex
cone $C \subset \R^4$.  This description can be extended to non-integer
parameters, which
then determine a size and shape for the hexagon, but not a triangulation.
Equation (\ref{eq: area of hexagon}) expresses the area,
measured in triangles, as a quadratic form of signature $(1,3)$. The
isometry group of any such a form is
$C_2 \times \Isom(\hy^3)$ (where $C_2$ denotes the
cyclic group of order $2$).

\begin{figure}[thbp]
\centering
\begin{minipage}{.55\textwidth}
	\includegraphics[width=\textwidth]{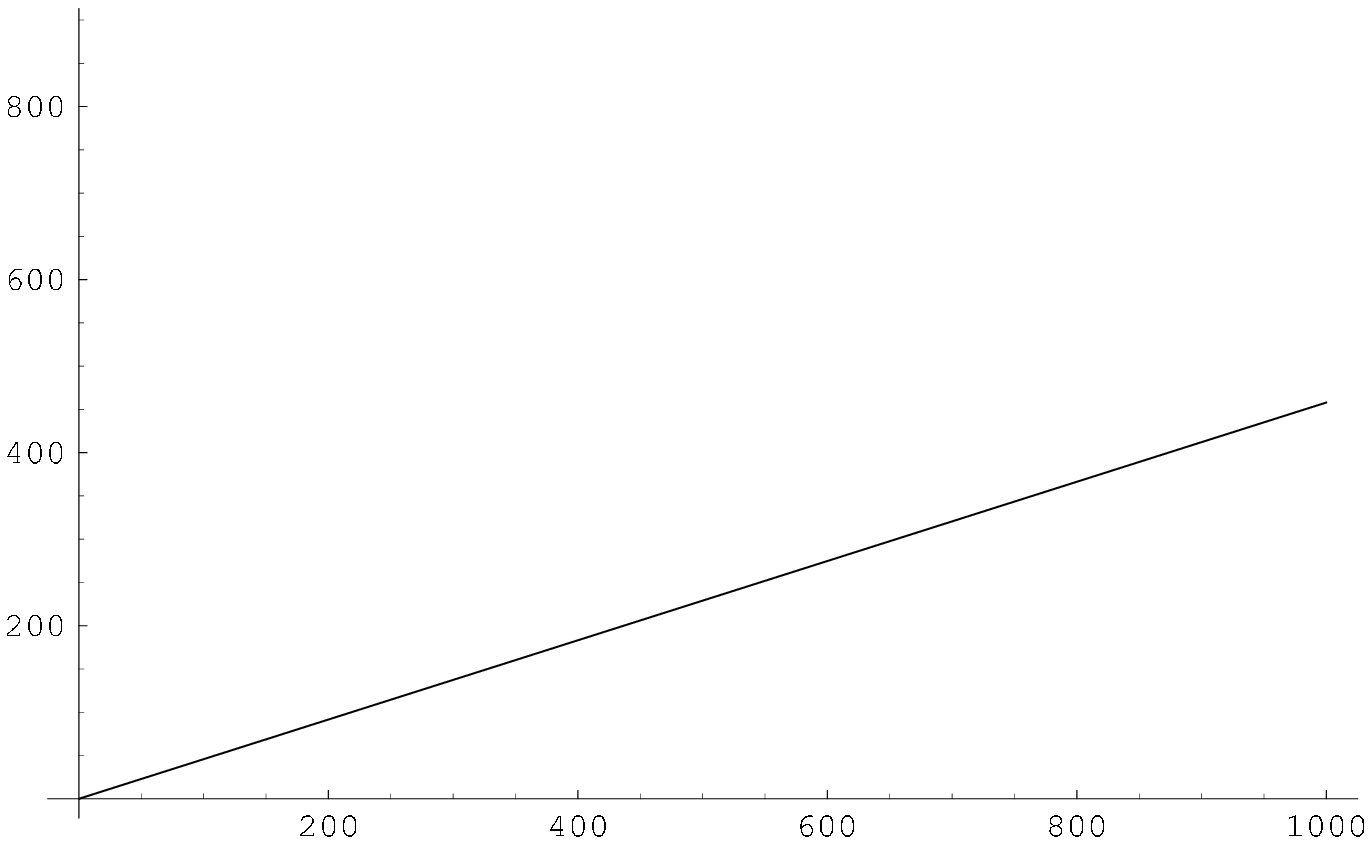}
\end{minipage}
\begin{minipage}{.44\textwidth}
\includegraphics[width=\textwidth]{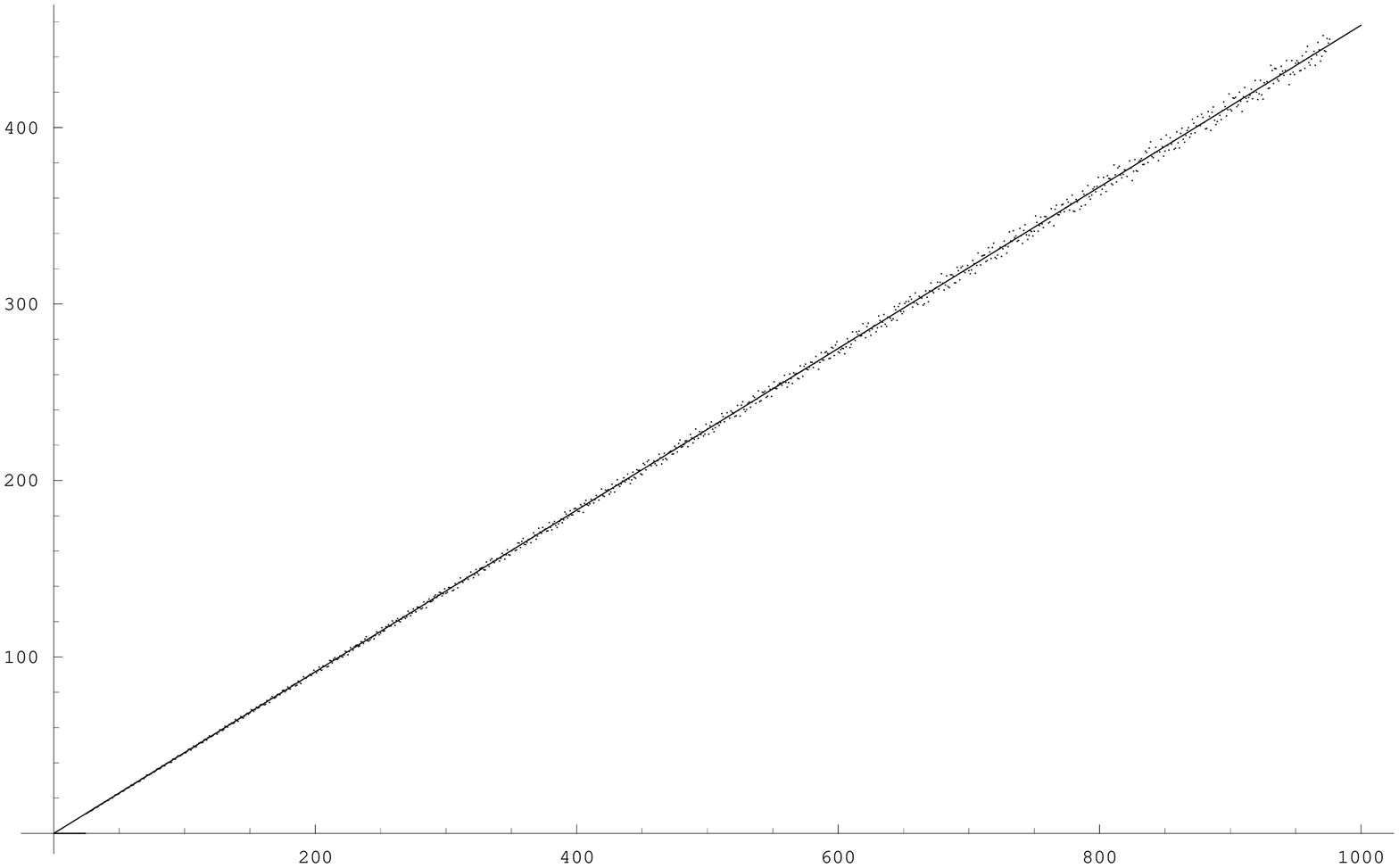}
\end{minipage}
\caption[Triangulations with up to 200 triangles.]
{\emph{Left}\qua the weighted count of Eisenstein lattice hexagons
containing up to 1000 triangles, using orbifold weights
$1/2^{k}$ where $k$ is the number of sides of a hexagon of length 0.
The parameter space of shapes (figure \ref{figure:
hexagon classification}) has hyperbolic volume 
$.91596559417\dots$ (1/4 that of the Whitehead Link complement), so
the number of hexagons containing $m$ triangles should grow 
on the average as the volume of the intersection of 
$C/2$ with the shell in $\euclidean^{(1,3)}$ between radius
$\sqrt{m}$ and $\sqrt{m+1}$, $.45798279709\dots * m$, as 
indicated. \emph{Right}\qua The same data averaged over windows 
of size 49.
\label{figure: hexagon enumeration}
}
\end{figure}
The possible shapes of lattice hexagons
(where rescaling is allowed) are parametrized by a convex polyhedron 
$H \subset \hy^3$ which is the projective image of the convex cone
$C \subset \reals^{(1,3)}$.
This polyhedron has three ideal vertices at infinity, which represent
the three directions in which shapes of hexagons can tend toward infinity,
by becoming long and skinny along one of three axes. In addition,
there are two finite vertices (top and bottom),
representing the two ways that a hexagon
can degenerate to an equilateral triangle. All dihedral angles of 
this hyperbolic polyhedron are $\pi/2$. Four edges meet at each
ideal vertex, while three edges meet at the finite vertices.
Triangulations with $m$ triangles are represented by a discrete
set $H_m \subset H$.  Figure \ref{figure: hexagon enumeration} plots
the count of how many of these lattice hexagons there are with each
possible area up to $1,000$. One indication of the relevance
of hyperbolic geometry is that the average number of hexagons of a given
area is well estimated by the hyperbolic volume of the parameter space.

\begin{figure}[htbp]
\centering
\begin{minipage}{.4\textwidth}
\centering
\includegraphics[width=.7\textwidth]{hex-triang.eps}
\end{minipage}
\begin{minipage}{.4\textwidth}
\centering
\includegraphics[width=.7\textwidth]{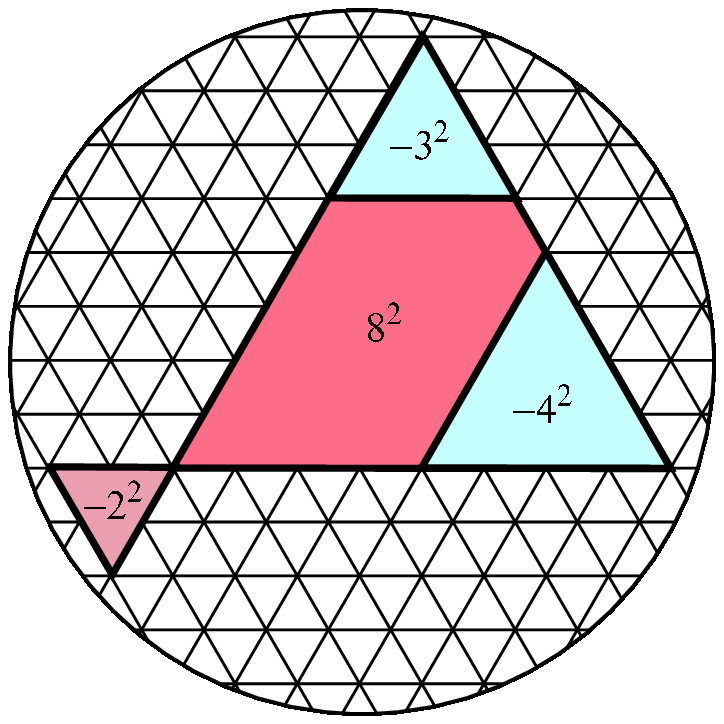}
\end{minipage}
\caption[Hexagon with butterfly.]{A butterfly operation moves
one edge of a hexagon (\emph{left}) across a butterfly-shaped quadrilateral
of 0 area, yielding a new hexagon (\emph{right}) of the same area.  The set
of butterfly moves generate a discrete group of isometries of
$\hy^3$, generated by reflections
in the faces of the polyhedron $H$.
}
\end{figure}
It's interesting to note
that $H$ is the fundamental polyhedron for
a discrete group of isometries of $\hy^3$, since all dihedral angles
equal $\pi/2$.  This group can be interpreted in terms of not necessarily
simple hexagons in the Eisenstein lattice  whose sides are parallel,
in order, to those of the standard hexagon. A non-simple lattice
hexagon wraps with integer degree around each triangle in the plane;
its total area, using these integer weights, is given by
the same quadratic form $n^2 - \sum p_i^2$.

Reflection in
a face of the polyhedron corresponds to a `butterfly move',
which is described numerically by reversing the sign of the length
of one of the edges of the hexagon, and adjusting the two
neighboring lengths so that the result is a closed curve.
Geometrically, the hexagon moves across a quadrilateral
reminiscent of a butterfly, resulting in
a new hexagon that algebraically encloses the same area as the original.
Note that this operation fixes any hexagon where the given side has 
degenerated to have length $0$---this is one of the faces of the
polyhedron $H$. The operations
for two sides of the hexagon that do not meet commute with each
other, and fix any shapes of hexagons where both these sides have
length $0$.  These shapes describe an edge of $H$, and since the
reflections in adjacent faces commute, the angle must be
$\pi/2$.  Two adjacent sides of the hexagon cannot both have
$0$ length at once, so the 9 non-adjacent pairs of sides of the
hexagon correspond 1--1 to the 9 edges of $H$.

Any solution to the equation $0 < m = n^2 - p_1^2 - p_2^2 - p_3^2$
determines a not necessarily simple hexagon of area $m$, which
projects to a point in $\hy^3$. By a sequence of
butterfly moves, this point can be transformed to be inside the
fundamental domain $H$.
The resulting point inside $H$ is uniquely determined by
the initial solution and does not depend on what sequence of butterfly
moves were used to get it there, since $H$ is
the quotient space (quotient orbifold) for the group action
as well as being its fundamental domain.

\section{Triangulations of the sphere}

Let $P(n; k_1, k_2, \dots, k_s)$ denote
the set of isomorphism classes of ``triangulations'' of
the sphere having exactly $2 n$ triangles, where for each $i$ there
is one vertex incident to $6 - k_i$ triangles, and all remaining
vertices are incident to $6$ triangles.
This paper will be limited to the non-negatively curved
cases that $ 0 < k_i \le 5$.
For there to be any actual
triangulations we must have $\sum_i k_i = 12$.
We will use the term ``triangulation'' throughout
to refer to a space obtained
by gluing together triangles by a pairing of their edges; thus,
in the case $k_i = 5$,
two edges of a triangle are folded together to form a vertex
incident to a single triangle.
Every triangulation of the sphere has an even number of triangles.

\begin{figure}[hbtp]
\begin{minipage}{.15\textwidth}
\centering
\includegraphics[width=.8\textwidth]{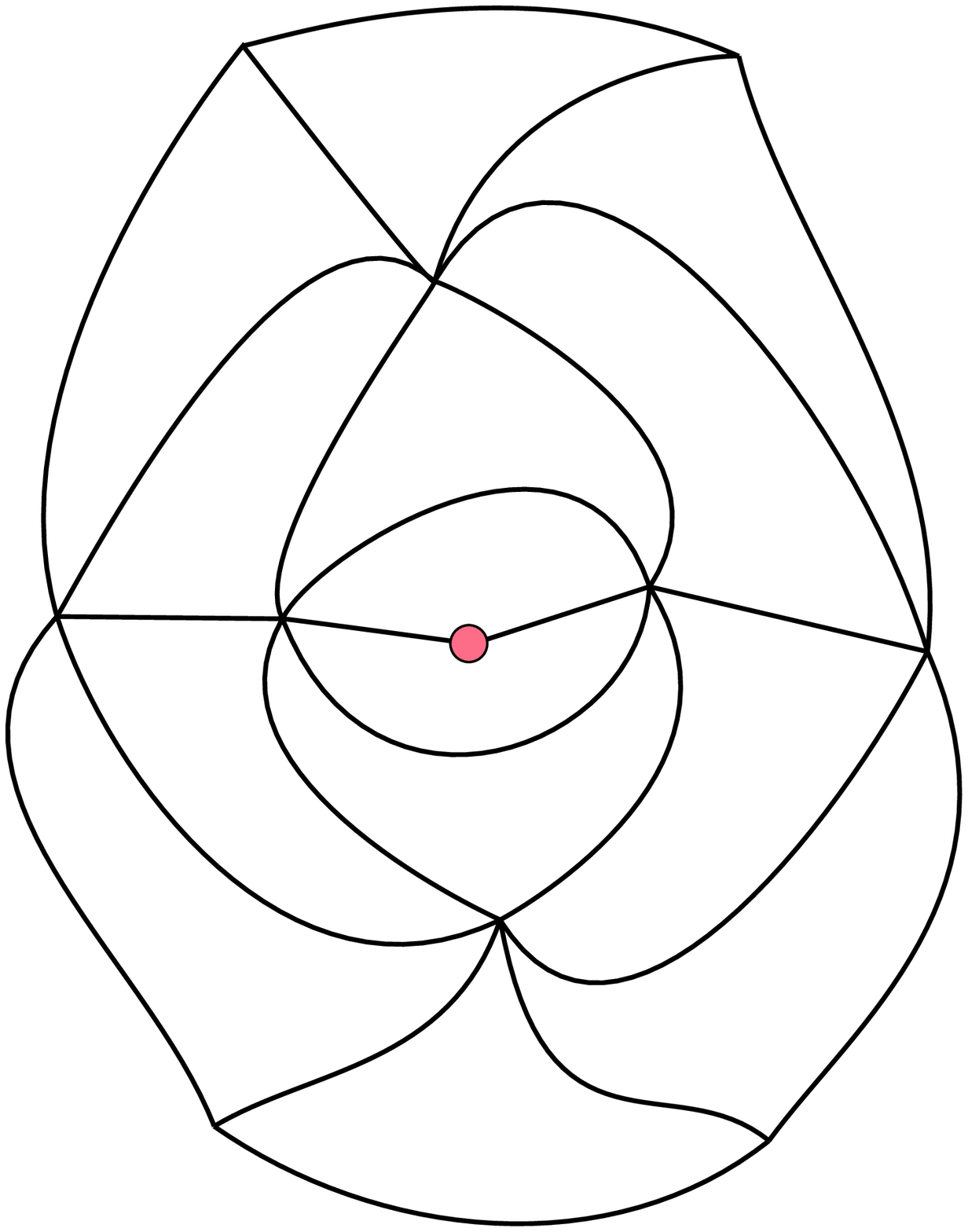}
\end{minipage}
\begin{minipage}{.68\textwidth}
\centering
\caption[The plane rolls three times around a vertex with two triangles.]
{A vertex with 2 triangles (\emph{left}) comes by folding up
a $120\degrees$ angle (\emph{right}). Equivalently, it has a neighborhood
whose developing map rolls around a vertex in the plane one time
for every  three revolutions of the cone.
Similar phenomena occur for vertices with 3 triangles or 1 triangle.
}
\end{minipage}
\begin{minipage}{.15\textwidth}
\centering
\includegraphics[width=.9\textwidth]{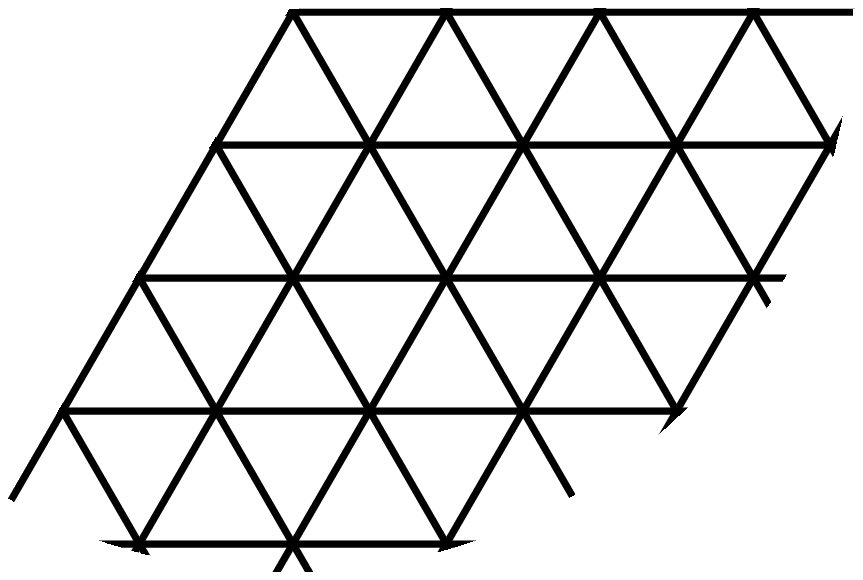}
\end{minipage}
\end{figure}

If $T \in P(n; k_1, \dots, k_s)$, then there is a
{\it developing map} $D_T$ from the universal cover $\tilde T$
of $T$ minus its singular
vertices into $E$.  Choose any triangle of $\tilde T$, and map it
to the triangle $\Delta ( 0, 1, \omega)$.  The developing map $D_T$ is
now determined by a form of analytic continuation, so that it
is a local isometry, mapping triangles to triangles.

\begin{floatingfigure}{1.1in}
\centering
\includegraphics[width=1in]{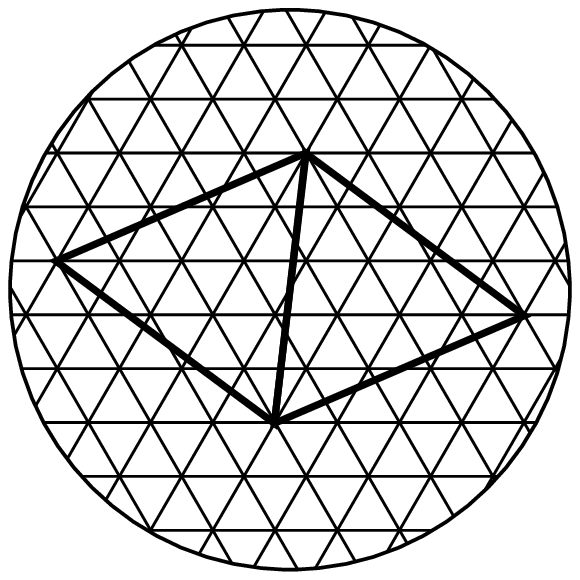}
\vskip .15in
\includegraphics[width=1in]{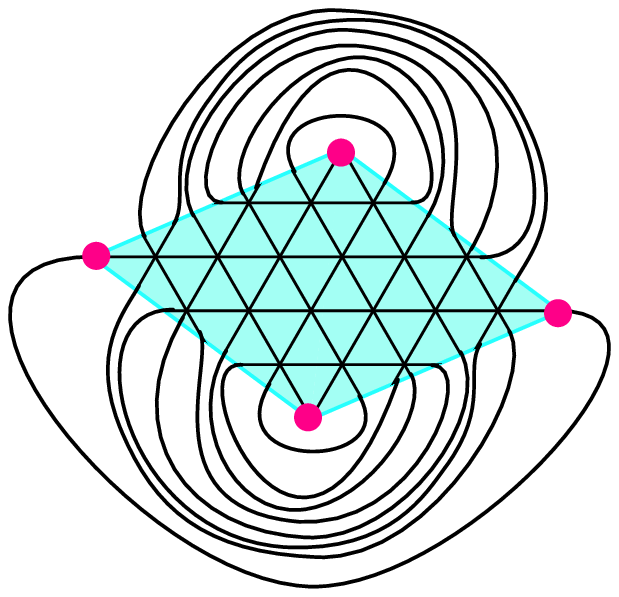}
\end{floatingfigure}
A particularly nice phenomenon happens for any vertices that have
only $1$, $2$, or $3$ triangles.
Consider a component $N_v$ of
the inverse image in $\tilde T$ of a small neighborhood of any such
vertex $v$.  It develops into the vicinity of some vertex $w$ in $\eisen$.
In these cases, the number of triangles around $v$ is a divisor
of $6$, so the developing map repeats itself when it first wraps around the
vertex $w$, along a path in $\tilde T$ which maps to a curve in $T$
wrapping respectively $6$, $3$, or $2$ times around the $v$.
Therefore, the developing map is defined from a smaller covering of
$T$ minus its singular vertices, which can be obtained as a certain
quotient space $S(T)$ of $\tilde T$.  In $S(T)$, each component
of the preimage of a small neighborhood of $v$ only intersects six
triangles.   In fact, $S(T)$ is isomorphic to $E$.   Therefore $T$ is
a quotient space of a discrete group $\Gamma(T)$
acting on $E$ such that only elements of $\eisen$ are
fixed points of elements of $\Gamma(T)$.

The examples where every vertex has 1, 2, 3 or 6 triangles are
$P(n;4,4,4)$, $P(n; 3,4,5)$ and $P(n; 3,3,3,3)$.
For $P(n; 4,4,4)$ or $P(n; 3,4,5)$, the group $\Gamma(T)$
is a triangle group.  A fundamental domain can be chosen as the union
of two equilateral triangles in the first case and
$30\degrees, 60\degrees, 90\degrees$ triangles of opposite orientation
in the second.  We may arrange that one of the vertices is at the origin.

Let $\alpha$ be a singular vertex closest to the origin.

In the case $T \in P(n; 4,4,4)$, 
the other singular vertices are $\eisen * \alpha$. Clearly this set
determines the group, and any $\alpha \ne 0$  will work.  The value
of $n$ is the ratio $\alpha \bar \alpha$
of a fundamental parallelogram
$0, \alpha, \alpha(1+\omega), \alpha\omega$
to the area of a primitive lattice parallelogram
$0, 1, 1+\omega, \omega$.  The possible numbers of triangles
are numbers expressible in the form $n = a^2 + 3 b^2$.

There is some ambiguity in this description:  if we replace
$\alpha$ by any of the other $5$ numbers $\omega^k \alpha$, we
obtain an isomorphic triangulation.   Thus, triangulations of this type
are in one-to-one correspondence with lattice points on the cone
$\complexes / \langle\omega\rangle$, where $\langle\omega\rangle$ 
refers to the multiplicative
subgroup of order $6$ generated by $\omega$.

Similarly, in the case $T \in P(n; 3,4,5)$, the vertices are of the form
$(m + p \sqrt {-3}) \alpha$, and  $n  = 2 \alpha \bar \alpha$.
As before, $\alpha$ is well-defined only up to multiplication by powers
of $\omega$.  In this case, if we replace $\alpha$ by $\omega^k \alpha$,
where $k$ is odd,
we get a different triangle group, but it has an isomorphic
quotient space.

\begin{figure}
\centering
\begin{minipage}{.3\textwidth}
\includegraphics[width=\textwidth]{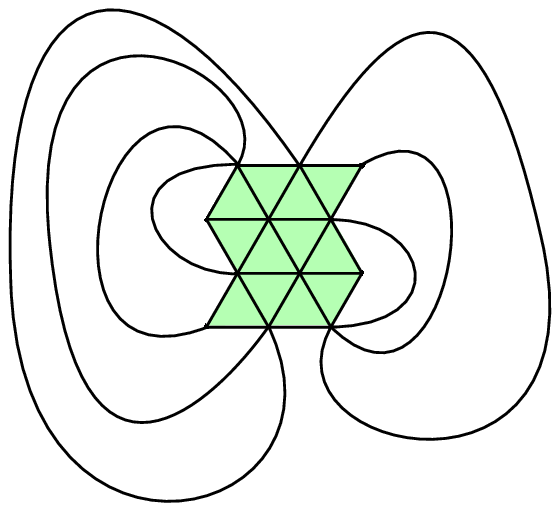}
\end{minipage}
\begin{minipage}{.3\textwidth}
\includegraphics[width=.8\textwidth]{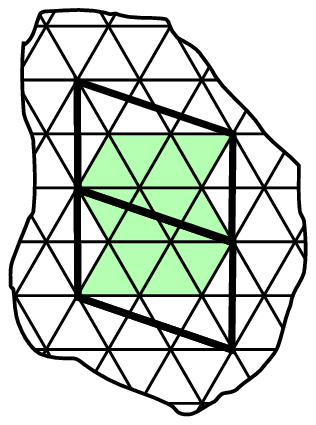}
\end{minipage}
\caption{Developing a triangulation with 3 or 6 triangles at
each vertex.}
\end{figure}

The case $P(n; 3,3,3,3)$ allows somewhat more variation.
For a singular vertex $x$ in $\eisen$, let $\gamma_x \in \Gamma(T)$
be the rotation of order $2$ about $x$.  Then for any two elements
$x$ and $y$, the product $\gamma_x \gamma_0 \gamma_y$ is a $180 \degrees$
rotation about $x + y$.  Therefore, the singular vertices form an
additive subgroup of $\eisen$.  Any additive subgroup will work.
The subgroup is determined if we specify the sides $\alpha$ and
$\beta$ of a fundamental parallelogram.  If we express $\alpha$
and $\beta$ as linear combinations of the generators $1$ and $\omega$
for $\eisen$, then the value of $n$ is twice the determinant of
the resulting two by two matrix. Every even number is achievable.
Of course, $\alpha$ and $\beta$ are well-defined only up to change
of basis for the lattice and up to multiplication by $6$th roots
of unity.  Note that multiplication by $\omega^3 = -1$ is also
represented by a change of generators.  A nice picture can be
formed by considering the shape parameter $z = \beta / \alpha$.
The action of the group $SL(2, \integers)$ on the set of shape
parameters is the usual action by fractional linear transformations
on the upper half plane.   Figure \ref{shapes of tetrahedra}
illustrates the set of shapes obtainable for $n = 246$.

Let us now skip to a more complicated case, that of
$P(n; 2,2,2,2,2,2)$, which includes the regular octahedron.
We have already encountered a special case: the hexagonal envelopes
of section \ref{triangulations} are examples of octahedra of this sort.

Just as a hexagon can be described by removing three small triangles
from a large triangle, there is a way to describe any element 
$T \in P(n; 2,2,2,2,2,2)$ by modifying an element
$\overline T \in P(m; 3,3,3)$, for some $m$.

Suppose $T$ is any triangulation of
the sphere with $6$ vertices incident to four triangles, and the rest
incident to $6$.

\cl{\hbox to 5pt{\hfil}
\includegraphics[width=.9\textwidth]{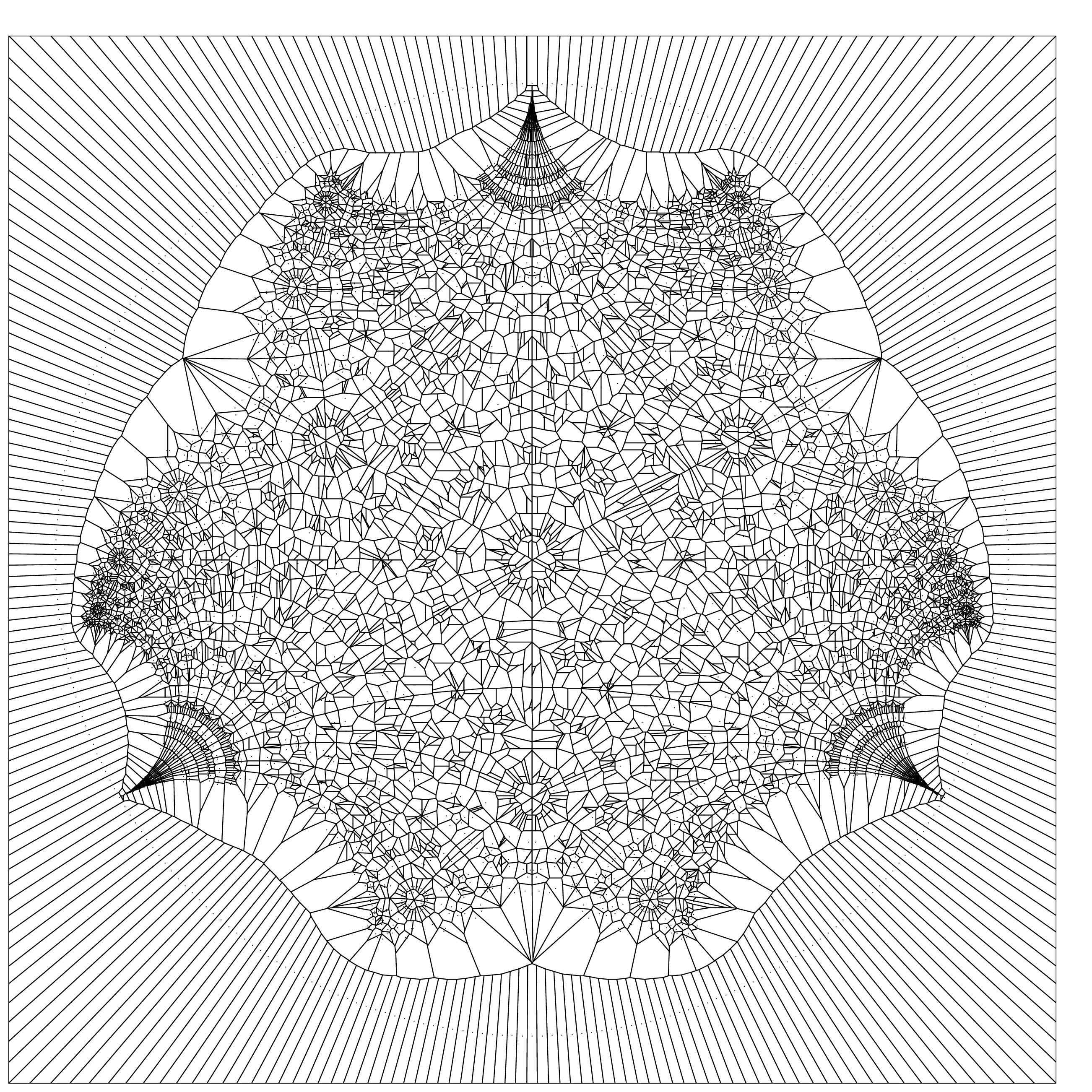}}\vglue -4mm
\begin{figure}[htbp]
\caption[Points in the hyperbolic plane describing triangulations of the sphere.]
{This is $P(246; 3,3,3,3)$,
plotted in the Poincar\'e disk model of $\hy^2$. The elements
of $P(246;3,3,3,3)$ are small dots; the Voronoi diagram for these
dots is shown, with one small dot inside each region. The position of
the dot in $\hy^2$ determines the
shape of a tetrahedron triangulated by 246 equilateral triangles.
Two dots which differ by $\PSL(2,\integers)$ represent the same shape.
The shape does not always completely determine the triangulation---one
also needs an angle for edges, that is, a lifting of the
point to a certain line bundle over $\hy^2$.
}
\label{shapes of tetrahedra}
\end{figure}

\begin{figure}[htbp]
	\begin{center}
	\begin{minipage}{.4\textwidth}
	\includegraphics[width=\textwidth]{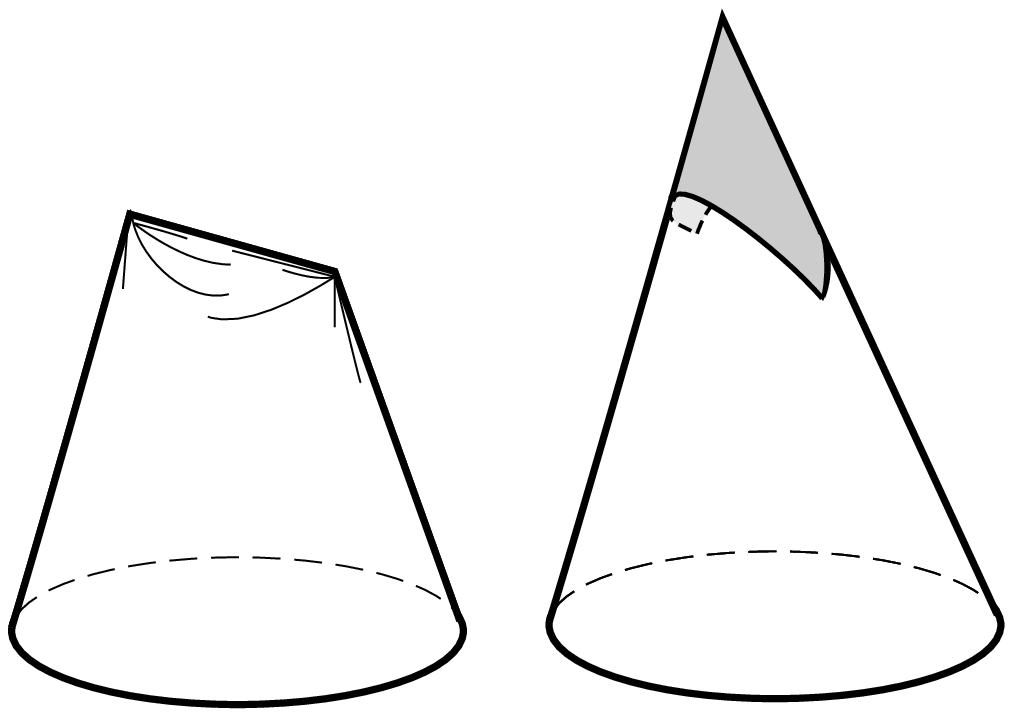}
	\end{minipage}
	\begin{minipage}
		{.15\textwidth}
	\end{minipage}
	\begin{minipage}{.35\textwidth}
		\includegraphics[width=\textwidth]{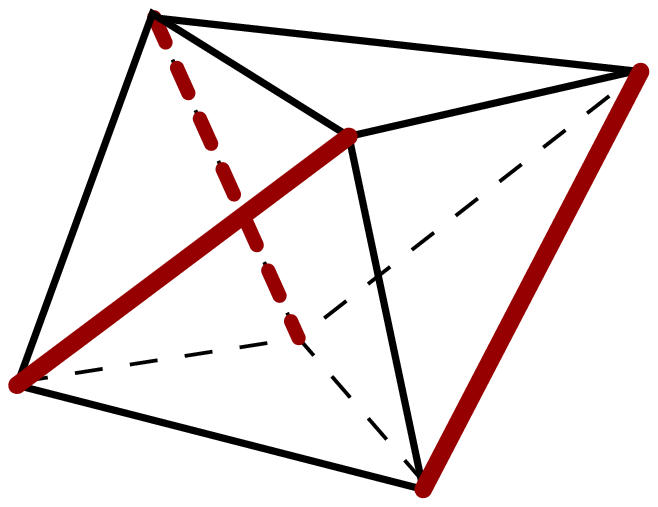}
	\end{minipage}
	\end{center}
\caption[A slit connecting two cone points can be filled in to
form a single cone.]{\emph{Left}\qua If a Euclidean cone manifold  is cut along
a geodesic arc joining the two cone points of curvature
$\alpha$ and $\beta$, the resulting figure is isometric to
a region in a Euclidean cone manifold with a new cone
point whose curvature is $\alpha + \beta$ (\emph{middle}). This gives
a recursive
procedure to reduce the construction of compact
Euclidean cone manifolds of non-positive curvature to ones
having only three cone points.
\emph{Right}\qua An element $T \in P(n; 2,2,2,2,2,2)$ can be 
reduced to $T' \in P(n'; 3,3,3)$ by slitting 3 arcs,
then extending. \label{figure: cutting a cone}
}

\end{figure}

Consider the associated cone metric $C$. 
We claim there is at least one way to join the $6$ cone points in pairs
by three disjoint geodesic segments.  To construct such a pairing,
first observe that any pair of cone points are joined by at least one geodesic:
the shortest path between them is a geodesic. Note that geodesics
can never pass through cone points with positive curvature,
except at their endpoints.  We see that there is
a collection of three not-necessarily disjoint geodesic segments joining
the $6$ points in pairs.  Let $\{ e, f, g\}$ be such a collection
of shortest possible length.  In particular, $e$, $f$ and $g$ are shortest
paths with their given endpoints.  No pair of these edges can intersect:
if they did, then by cutting and pasting, one would find that the four
endpoints involved could be joined in an alternate way by shorter paths.

\begin{figure}[ht]
\centering
\includegraphics[width=.6\textwidth]{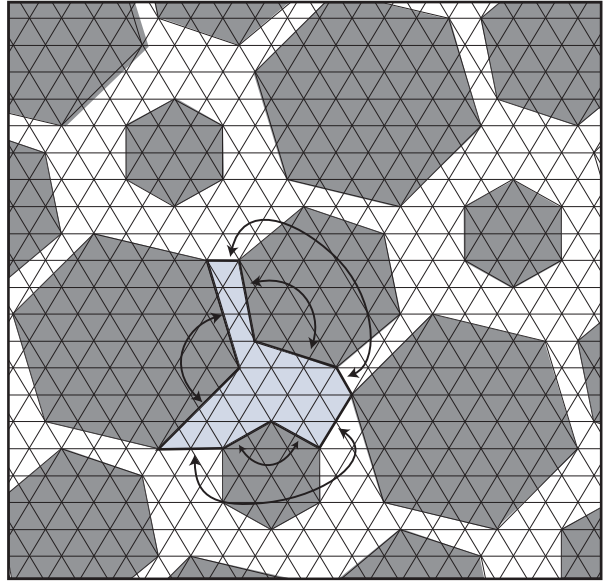}
\caption[A cut-out diagram for a generalized octahedron.]{This is an illustration of the construction of a generalized octahedron,
that is, an element of $P(n; 2,2,2,2,2,2)$.  First,
choose a $3,3,3$ group acting
in the plane with the fixed points of the elements
of order $3$ on lattice points of $\eisen$.
Then choose three families of lattice hexagons invariant by the group,
centered at the fixed points of elements of order $3$.  Remove the
hexagons, form the quotient by the group, and glue the edges of the
resulting slits together. Equivalently, you can glue the boundary
of a fundamental domain as illustrated. }
\label{constructing octahedra}
\end{figure}

Cut $C$ along the three edges $e$, $f$ and $g$, and consider the
developing map for the resulting surface $C'$.  At an endpoint of say
$e$, the developing image subtends an angle of $120 \degrees$; a curve
which wraps three times around $e$ in a small neighborhood
develops to a curve wrapping once around 
the outside of a regular hexagon $H_e$ in the
plane.  Let $C_e$ be $H_e$ modulo a rotation of order $3$.  If we glue
$C_e$ and the similarly constructed cones $C_f$ and $C_g$ to the cuts,
we obtain a new cone-manifold $C''$, with three cone points of order $3$.
The hexagon $H_e$ has its vertices on lattice points of $\eisen$, so its
center is also a lattice point of $\eisen$.  Therefore, $C'' \in P(m; 4,4,4)$
for some $m$.  Consequently,
a general element of $P(n; 2,2,2,2,2,2)$ is
obtained by choosing some $m$ bigger than $n$, choosing an element of
$P(m; 4,4,4)$, and choosing three types of hexagons whose area in triangles
adds to $6(m - n)$ such that when they are placed around the three classes
of order $3$ points in the plane, all their images are disjoint.  Cut all
these hexagons out of the plane, divide by the $(3,3,3)$ triangle group, and
glue together the pair of edges coming from each hexagon.
We can express this as a choice of four elements $\alpha_i \in \eisen$,
such that 
\[
\alpha_1 \bar \alpha_1 - \alpha_2 \bar \alpha_2 - \alpha_3 \bar \alpha_3 -
\alpha_4 \bar \alpha_4
: \]
$\alpha_1$ is used to construct
the original $(3,3,3)$ triangle group, and the other $\alpha_i$'s are
vectors from the centers of the each of the hexagons to one of the
vertices, yielding a triangulation 
$$ T(\alpha_1, \alpha_2, \alpha_3, \alpha_4) \in P(n; 2,2,2,2,2,2).$$
The $\alpha_i$'s are
subject to an additional geometric condition, that the hexagons they
define be embedded.  The coordinates are only defined up to a
geometrically-defined equivalence relation, having to do with the
multiplicity of choices for $e$, $f$, and $g$.
The easy observation is
that when any of the $\alpha_i$ are multiplied by powers of $\omega$,
we obtain the same $T$.  These coordinates make it easy to automatically
enumerate all examples, although it is somewhat harder to weed out
repetitions.  The geometric conditions can be nearly determined
from the norms: if  $|\alpha_i| + |\alpha_j| < |\alpha_1|$,
for $i \ne j \in \{2,3,4\}$, then the hexagons are disjoint; if this sum
is greater than $(2 / \sqrt{3}) |\alpha_1| = 1.1547\dots |\alpha_1|$,
then two hexagons intersect; otherwise, one needs to consider the picture.
If $|\alpha_i| < |\alpha_1|/3$ for $i > 1$, then the three edges $e$, $f$
and $g$ are clearly the three shortest possible edges; in general, the
question is more complicated.
The standard octahedron $O \in P(4; 2,2,2,2,2,2)$, for example, has
an infinite number of
descriptions, for example
$O = T(2 k + 1 + (-k+2) \omega, k+\omega, k+\omega, k+\omega)$
for every $k \ge 0$.

Another construction will be given in section
\ref{section: construction and fundamental domain} that can be used
to search all possibilities while weeding out repetition fairly
efficiently.

\section {Shapes of polyhedra} 

Any collection of $n$--dimensional Euclidean
polyhedra whose $(n-1)$--dimensional faces are glued together
isometrically in pairs yields
an example of a {\it cone-manifold}
and gives a pretty good flavor for
the singular behavior that can occur.
However, polyhedra are
not a suitable substrate for a definition in the context we need,
since we will be working with metrics whose local geometry
has no concept of polyhedra comparable to the Euclidean case:
they have no totally geodesic hypersurfaces.

In general, a cone-manifold is a kind of singular Riemannian metric;
in our case, we will work with
spaces modelled after a complete Riemannian $n$--manifold
$X$ together with a group $G$ of isometries of $X$,
called an $(X, G)$--manifold.  If $G$ acts transitively, this would be
called a \df{homogeneous space}, but {\it 
$G$ does not necessarily act transitively}. Moreover, the
group $G$ is part of the structure. It is not necessarily
the full group of isometries of $X$: for instance, we might have
$X = \euclidean^2$ and $G$ the group of isometries that preserve
the $\eisen$.

 An {\it $(X, G)$--manifold} is a space
equipped with a covering by open sets with homeomorphisms
into $X$, such that the transition maps on the overlap of any
two sets is in $G$.

The concept of an $(X,G)$--cone-manifold is defined inductively by
dimension, as follows:

If $X$ is $1$--dimensional, an $(X,G)$--cone-manifold
is just an $(X,G)$--manifold.

Suppose $X$ is $k$--dimensional, where $k > 1$.  For any point $p \in
X$, let $G_p$ be the stabilizer of $p$, and let $X_p$ be the set of
tangent rays through $p$. Then $(X_p, G_p)$ is a model space of one
lower dimension.  If $Y$ is any $(X_p, G_p)$--cone-manifold, there is
associated to it a fairly intuitive construction, the \df{radius $r$
cone} of $Y$, $C_r(Y)$ for any $r > 0$ such that the exponential map
at $p$ is an embedding on the ball of radius $r$ in $T_p(X)$,
constructed from the geodesic rays from $p$ in $X$ assembled in the
same way that $Y$ is.  That is, for each subset of $X_p$, there is
associated a cone in the tangent space at $p$, and to this is
associated (via the exponential map) its radius $r$ cone in $X$.
These are glued together, using local coordinates in $Y$, to form
$C_r(Y)$.

An $(X,G)$--cone-manifold is a space such that  each point
has a neighborhood modelled on the
cone of a compact, connected $(X_p, G_p)$--manifold. 

One reason for considering inhomogeneous model spaces $(X,G)$
is  that even if we start with an example as homogeneous as
$(\CP^n, U(n))$, during the inductive examination of tangent
cones we soon encounter 
model spaces $(X,G)$ where $G$ is not transitive.  

If $C$ is an $n$--dimensional
$(X,G)$--cone-manifold, then a point $p \in C$ is a {\it regular}
point if $p$ has a neighborhood equivalent as an $(X,G)$--space to a
neighborhood in $X$, otherwise it is {\it singular}.  It follows
by induction that regular points are dense, and that $C$ is the metric
completion of its set of regular points.
The distinction between regular points and singular points can be refined
to give the concept of the {\it codimension} of a point $p \in C$.
If the only cone type neighborhood that a point $p$ belongs to
is the neighborhood
centered at $p$, then $p$ has codimension $n$.  Otherwise, there is
some cone neighborhood centered at a different point $q$ that $p$
belongs to, and the codimension
of $p$ is defined inductively to be the codimension of the ray through
$p$ in $(X_q,G_q)$.

By induction, it follows that
each point $p$ of codimension $k$ is on an $(n-k)$--dimensional
stratum of $C$ which is locally isometric to a totally
geodesic subspace $E_p \subset X$ ---
this stratum is an $(E_p, G(E_p))$--space, where $G(E_p)$ is the subgroup
of $G$ sending $E_p$ to itself.

An oriented Euclidean, hyperbolic, or elliptic cone-manifold of dimension $n$ is
a space obtained from a collection of totally geodesic simplices
via a $2$ to $1$ isometric
identification of their faces.

Suppose that $n$ numbers $\alpha_i$ are specified,
all less than $1$, such that $\sum \alpha_i = 2$.  Let
$C(\alpha_1, \alpha_2, \dots, \alpha_n)$ be the space of Euclidean
cone-manifold structures on the sphere with
$n$ cone singularities of curvature $\alpha_i$ (cone
angles $2 \pi  ( 1 - \alpha_i) $), up to equivalence by
orientation-preserving similarity.
We do not specify any homotopy class of map relative to the cone points,
nor any labelling of the cone-points in these equivalences.  Let
$P(A; \alpha_1, \dots, \alpha_n)$ be the finite-sheeted covering in
which the cone points can be consistently labelled.  Note that the
fundamental group of $P(A; \alpha_1, \dots, \alpha_n)$ is the
pure braid group of the sphere, and the fundamental group of $C(\alpha_1,
\dots, \alpha_n)$ is contained in the full braid group of the sphere and
contains the pure braid group.  The exact group depends on the collection
of angles, since only cone points with equal angles can be interchanged.

How can we understand these spaces?  We will first
construct a local coordinate system for the space of shapes of such
cone-metrics, in a neighborhood of a given metric $g$.

\begin{proposition}[Cone-metrics have triangulations]\label{cone-metrics have triangulations}
Let $C$ be any metric on the sphere which is locally Euclidean except
at isolated cone-points of positive curvature.  Then $C$ admits a
triangulation in the sense of a subdivision of $C$ by images of
geodesic Euclidean triangles, possibly with identifications of
vertices and/or edges, with vertex set the set of cone points.
\end{proposition}
\begin{proof}  Associated to each cone point
$v$ of $C$ is the open \emph{Voronoi region}
for $v$, consisting of those points $x \in C$ which are closer 
to $v$ than to any other cone point, and furthermore, 
have a unique shortest geodesic arc connecting $x$ to 
$v$.  A 
\emph {Voronoi edge} consists of points $x$ that have 
exactly two shortest geodesic arcs to cone points.  Each 
Voronoi edge is a geodesic segment.  
It can happen that a Voronoi edge has the same Voronoi region on
both sides
if $C$ has a fairly long, skinny region with a cone point $v$ 
far from other cone points. Take any point $x$ on a Voronoi edge,
and let $D$ be the largest metric ball centered at $x$ whose interior 
contains no cone points. Then $D$ is the image of an isometric immersion
of a Euclidean disk $D'$ with exactly
two points $v_1, v_2 \in \partial D'$ that map to cone points of $C$.
The chord $\overline{v_1 v_2}$ of $D'$ maps to an
arc in $C$. The collection of all such arcs have disjoint interiors,
for if not, one could lift the situation to $\euclidean^2$: 
whenever two chords of two distinct disks in $\euclidean^2$ cross,
at least one of the four endpoints is
in the interior of at least one of the two disks.

The \emph {Voronoi vertices}
are those points that have 
three or more shortest arcs to cone points. The largest metric disk
about a  Voronoi vertex with no cone points in the interior is the 
image of an isometrically immersed Euclidean disk.
The convex hull of the set of points on the 
boundary of the Euclidean disk that map to cone points 
is a convex polygon mapping to $C$ with boundary mapping to 
the edges previously constructed.   Subdivide each of these 
polygons into triangles by adjoining diagonals.  The result 
is a geodesic triangulation of $C$ in the sense of the 
proposition whose vertex set is the 
set of cone 
points.
\end{proof}

Let $T$ be any geodesic triangulation of the cone-manifold $C$;
it might or might
not be obtained by this construction.  Choose one of the edges of $T$,
and map it isometrically into $\complexes$, with one endpoint at the origin.
This map extends to an isometric developing map $D\co  \tilde C \arrow \complexes$,
where $\tilde C$ is the universal cover of the complement $C_0$ of the vertices
of $C$.  Associated with
each directed edge $e$ of the triangulation $\tilde T$ of
$\tilde C$ is
a complex number $Z(e)$ (really a vector), the difference between its endpoints.
These vectors satisfy the cocycle condition, that the sum of the vectors
associated to the oriented boundary of a triangle is $0$.
Let $H\co  \fund {C_0} \arrow \isom ( \euclidean^2)$
be the holonomy of the Euclidean structure, and let
$H_0\co  \fund {C_0} \arrow S^1 \subset \complexes$
be its orthogonal part.  If $\tau_\gamma$ is the covering transformation
of $\tilde C$ over $C_0$ associated with the element $\gamma \in \fund {C_0}$,
then $Z(\tau_\gamma (e)) = H_0 ( \gamma) Z(e)$.  In other words, it is
a cocycle with twisted coefficients --- the coefficient bundle
is the tangent space of $C_0$.   Euclidean structures near $C$,
up to scaling, are
parametrized by cocycles near $Z$, up to multiplicative complex numbers,
since any nearby cocycle determines a collection of shapes of triangles
which can be glued together to form a cone-manifold with the same
set of cone angles.

It is clear that change of coordinates, from those given by $T$
to those given by a triangulation $T'$, is a linear map, since the
developing map for the edges of $T'$ can be computed as a linear function
of a cocycle expressed in terms of $T$.

\begin{proposition}[Dimension is $n-2$]\label{dimension is n-2}
The complex dimension of the space of cocycles, as described above, is
$n - 2$, where $n$ is the number of vertices.
\end{proposition}

See \cite{MR96h:57010} for various computations related 
to this.

\begin{proof}
We will describe a concrete construction for a basis
for the cocycles, which amounts to making a gluing diagram to construct
$C$ from a polygonal region on a cone.\footnote{In the general, complicated
cases, this would likely be an immersed polygonal region on a cone.}

We will divide the set of edges into {\it leaders} (the basis elements) and
{\it followers}.
Begin by picking any vertex $v_{last}$ of $T$,
and designate all edges leading into
that vertex as followers.
Now pick a tree in the $1$--skeleton connecting all vertices
except $v_{last}$: these will be leaders. The remaining edges
are additional followers.
There is a dual tree, in the dual $1$--skeleton of the cell-division formed
by removing the followers touching $v_{last}$, consisting of the $2$--cells
and the remaining followers.

Suppose the value of a $1$--cocycle is specified on each of the leaders.
We can then calculate it on each of the followers, as follows.  Inductively,
if the current
dual tree of undetermined values is bigger than a single point, pick
a leaf of the tree.  This is a follower which is part of a triangle whose
other two sides have determined values; from them,
we determine the value for the follower to satisfy the coboundary condition
on the given triangle.  What remains is still a tree.

Finally, we are left with everything determined, except for $v_{last}$
and its remnant cluster of followers.  At this point, we have enough information
to determine the affine holonomy around $v_{last}$.  The orthogonal
part is a non-trivial rotation, so that
it has a unique fixed point.  The values of the cocyle for the remaining
followers are determined by pointing them toward the fixed point.

A spanning
tree for the $n-1$ vertices excluding the last has $n-2$ edges, so
the space of cocycles is $\complexes ^ {n-2}$. The projective space then has
dimension ${n-3}$.
\end{proof}

The area of a cone-manifold structure defines a hermitian form on the
space of cocyles: that is, given a cocycle $Z$,
$A(Z) =  \frac{1}{4} \sum_{\rm triangles} i e_1 \bar e_2 - i e_2 \bar e_1$
where in local coordinates $e_1$ and $e_2$ are successive edges of the triangle
proceeding counterclockwise.  Obviously $A(Z)$ is independent of
choice of local coordinates.

\begin{proposition}[Signature $(1,n-3)$]\label{signature(1,n-3)}  If each of the
$\alpha_i > 0$, then $A$ is a hermitian form of signature
$(1,n-3)$.
\end{proposition}
\begin{proof}
We have seen this illustrated in several examples already.  There is
a general procedure for diagonalizing the expression for area.
If $C$ has only three vertices, then the vector space is only
one dimensional, so $A$ is necessarily positive definite:
it is proportional to the square of the length of any of the edges
of $T$.

We have already seen
the special case that there are four cone angles all equal to
$\pi$, under the guise of $P(n; 3,3,3,3)$.
The expression for area is the determinant of
a $2 \cross 2$ real matrix, made of the real and imaginary parts
of two of the values $Z(e)$.
Since determinants
can be positive or negative, this is a hermitian form of signature $(1,1)$.

In every other case, there are at least two cone angles whose curvatures
have sum less
than $2 \pi$.  Construct any geodesic path $e$ between them, slit $C$ open, and
glue a portion of a cone with curvature the sum of the two curvatures to
obtain a cone-manifold $C'$ with one fewer singular points
(figure \ref{figure: cutting a cone}).
The area of $C$ is the area of $C'$ minus a constant
times the square
of the length of $e$.  This gives an inductive procedure for diagonalizing
$A$, inductively showing that the signature of the
area is $(1,n-3)$.
\end{proof}

The set of positive vectors in a Hermitian form of signature $(1,n-3)$
up to multiplication by scalars, is biholomorphic to the interior of
the unit ball in
$\complexes ^ {n-3}$, and is known as {\it complex hyperbolic space}
$\CH^{n-3}$.
A metric of negative curvature is induced from the Hermitian form;
as a Riemannian metric, its sectional curvatures are pinched between
$-4$ and $-1$.  Therefore, $C (A; \alpha_1, \dots, \alpha_n)$
is a complex hyperbolic manifold.

It is not metrically complete, however.
Any two singular points of a  $c$
whose curvature adds to less than $2 \pi$ 
can collide as the cone-metric changes a finite amount, measured in
the complex hyperbolic metric.  We will next examine how to adjoin to
$C(\alpha_1, \dots, \alpha_n)$ the degenerate cases where one
or more of the cone points collide, to obtain a space
$\bar C(\alpha_1, \dots, \alpha_n)$ which is the metric completion
of $C(\alpha_1, \dots, \alpha_n)$.

Each element $c$ of $\bar C(\alpha_1, \dots, \alpha_n)$ is associated with
some partition $P$ of the angles $\alpha_i$;  $c$ is a Euclidean
cone-manifold where each cone point is associated with
a partition element $p \in P$ and has curvature equal to the sum of
the elements of $p$.  We regard two partitions as equivalent if one
can be transformed to the other by a permutation of the index set
which preserves the values of the $\alpha_i$.
A limit of
a sequence of cone-manifolds associated with some partition will be
associated with a coarser partition, if distances between
some of the cone points in the sequence tend to zero.

\begin{theorem}[Completion is cone-manifold]\label{Completion is cone-manifold} 
The metric completion of\break $C(A, \alpha_1, \dots, \alpha_n)$
is $\bar C(\alpha_1 , \dots, \alpha_n)$, which
is a complex hyperbolic cone-manifold.
\end{theorem}

\begin{proof}  
There is a very natural way to describe
regular neighborhoods for the stratum $S_P$ corresponding to a partition
$P$ of the set of curvatures concentrated at cone points.

Consider an element $c \in C(\alpha_1, \dots, \alpha_n)$
such that the cone points are
clustered in accordance with $P$.  We may assume that the diameter
of each cluster is less than the minimum distance from the cluster
to any cone point not in the cluster, and less than some small
constant $\epsilon$.

The holonomy
for a curve which goes around any cluster $D$ is a rotation by the total
curvature of $D$, unless the total curvature is $2 \pi$. When the total
curvature of $D$ is $2\pi$, the holonomy is a translation.
If the holonomy is actually a rotation,
it leaves invariant each of a family of circles; with our assumption
that the cluster is isolated from other cone points, the encircling
curve is isotopic to one of these circles.

If the total curvature of $K(D)$ is less than
$2\pi$, the surface of $c$ near such a circle isometrically
matches a cone with apex on
the same side as the cluster, with cone point of curvature $K(S)$.
In this case, we can define a new cone-manifold
$p(c)$ by cutting out each cluster, and replacing it by a portion of this cone.
In local coordinates, this gives a local orthogonal
projection from a neighborhood of $c$ to $S_P$.  The distance
from the singular stratum is $\sqrt{ \area p(c) - \area c}$.
Note that the normal fibers for strata corresponding to
 subclusters of a cluster are contained in
normal fibers for the larger cluster.

\begin{figure}[htbp]
\begin{minipage}{.4\textwidth}
\includegraphics[width=\textwidth]{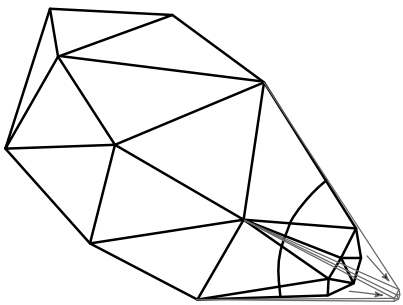}
\end{minipage}
\begin{minipage}{.05\textwidth}\hfil
\end{minipage}
\begin{minipage}{.53\textwidth}
\caption[A cluster of cone points shrinks to one cone point.]
{Any cluster of cone points close together compared
to the distance to other cone points can be shrunk to a single cone 
point. This process gives a radial structure to a neighborhood of a
singular point in the space of cone-metrics with designated curvatures
on a sphere.}
\end{minipage}
\end{figure}
The total curvature cannot be greater than
$2\pi$, if $\epsilon$ is chosen properly:
in that case, $c$ would match
the surface of a cone with apex on the opposite side from the cluster.
The area of $C$ is less than the area of the portion of cone, plus the area
within the cluster, so that
if $\epsilon$ is small compared to $\theta/A$, where $\theta$ is the
minimum value by which a curvature sum can exceed $2 \pi$, 
this cannot occur.

A cluster of arbitrarily small diameter with total curvature $2 \pi$ can occur,
but this forces the diameter of $c$ to be large: in this case, $c$ matches
the surface of a cylinder outside a neighborhood of the cluster, and there
is a complementary cluster at the other end of the cylinder.
As $c$ moves a finite distance in the complex hyperbolic metric,
its diameter cannot go to infinity, so no such cluster goes to $0$ in diameter
in the metric completion of $C(\alpha_1, \dots, \alpha_n)$.

Within any bounded set of $C(\alpha_1, \dots, \alpha_n)$,
we are left only with the case of small diameter
clusters whose total curvature is less than $2\pi$.

It is now easy to see that $\bar C(\alpha_1, \dots, \alpha_n)$
is the metric completion of $C(\alpha_1, \dots, \alpha_n)$ and that
it is a complex hyperbolic cone-manifold.
\end{proof}

Of particular importance are the strata of complex codimension $1$ or
real codimension $2$.  These strata correspond to the cases when only
two cone points of $c$ have collided.  What are the cone angles around
these strata?

\begin{proposition}[Cone angles around collisions]\label{cone angles around collisions}
Let $S$ be a stratum of\break $C(\alpha_1, \dots, \alpha_n)$ where two
cone points with curvature $\alpha_i$ and $\alpha_j$ collide.

If $\alpha_i = \alpha_j$, the cone angle $\gamma(S)$ around $S$ is
$\pi - \alpha_i$, otherwise it is
$2\pi - \alpha_i - \alpha_j$.

\end{proposition}

In other words, the cone angle in parameter space
is the same as the physical angle two nearby cone points go through,
as measured from the apex of the cone that would be formed by their collapse,
when they revolve about each other
until they return to their original arrangement.

\begin{proof}
When cone points $x_i$ and $x_j$ with these two angles are close together on
a cone manifold $c$,
we can think of $c$ as constructed from $p(c)$ by replacing a small
neighborhood of the cone by a portion of a modified cone
$D(\alpha_i, \alpha_j)$
with two cone points.  The shape of $D(\alpha_i, \alpha_j)$
is uniquely determined by $\alpha_i$ and $\alpha_j$
up to similarity.  Thus, the shape of $c$ is determined by selecting
the point $x_i$ on $p(c)$, and may be represented by $p(c)$ together
with the vector $V$ from the combined cone point of $p(c)$ to $x_i$.
In local inhomogeneous
coordinates coming from a choice of a triangulation, $V$ is a locally
linear function, described by a single complex number.

If $\alpha_i = \alpha_j$, then when the argument of $V$ is increased
by half the cone angle, or $\pi - \alpha_i$, $x_i$ and $x_j$ are interchanged,
and the resulting configuration is indistinguishable.  Therefore,
$\pi - \alpha_i$ is the cone angle along $S$, (and $\pi + \alpha_i$ is
the curvature concentrated at $S$).
If $\alpha_i \ne \alpha_j$,
the argument of $V$ must be increased by the cone angle, $2 \pi - \alpha_i -
\alpha_j$, before the same configuration is obtained again.  In this case,
$2 \pi - \alpha_i - \alpha_j$ is the cone angle along $S$,
and $\alpha_i + \alpha_j$ is the curvature concentrated along $S$.
\end{proof}

More generally, if $S$ is a stratum of complex codimension $j$ representing
the collapse of a cluster of $j+1$ cone points, each normal fiber to $S$
is a union of `complex rays', swept out by an ordinary real ray
by rotating it the direction $i$ times the radial direction. The 
space of complex rays is the complex link of the stratum,
a complex cone-manifold whose complex dimension
is one lower. The real link is a Seifert fiber space over
the complex link, with generic fiber a circle of length $\alpha$
which we can call the \emph{scalar cone angle} $\gamma(S)$ at $S$.
We define the \emph{real link fraction} of $S$ to be the ratio of the
volume of the real link of $S$ to the volume of $S^{2j-1}$
(the real link in the non-singular case), and
similarly the \emph{complex link fraction} is the ratio of the
volume of the complex link to the volume of $\CP^{j-1}$.

\begin{proposition}[Cone angles for multi-collisions]
\label{proposition: multiple cone angles}
Let $S$ be a stratum of complex codimension $j$
where $j+1$ cone points of curvature $\kappa_1, \dots, \kappa_j$
collapse.

Let $N$ be the order of the subgroup of the symmetric group
$S_j$ that preserves these numbers.  Then:
\begin{itemize}
\item[\bf a\rm)]
The scalar cone angle is
\[ \gamma(S) =  2 \pi - \sum_i \kappa_i . \]
\item[\bf b\rm)]
The complex link fraction is
\[ \frac {(\gamma(S)/2\pi)^{j-1} }{N} . \]
\item[\bf c\rm)]
The real link fraction is
\[ \frac {(\gamma(S)/2\pi)^j}{N} .  \]
\end{itemize}
\end{proposition}
\begin{proof}
The proof of part (a)
is the same as above, with the observation that a cluster of 3 or more
cone points can always be slightly perturbed to make it asymmetrical,
so in the generic fiber of the Seifert fibration (obtained by rotating
the cluster of cone points) no permutations
of the cone points occur.

For (b), think first about the case that all cone angles are different,
so as to avoid a symmetry group at first.  A neighborhood
of $S$ is then a manifold, isomorphic to the limiting case when
$\kappa_i \to 0$, the space
of $(j+1)$--tuples in the plane up to affine transformations. The
complex link is a complex cone-manifold structure on $\CP^{j-1}$.
If $\omega$ is a closed $2$--form on $\CP^{j-1}$ that integrates to
$1$ over $\CP^1$, then $\omega^{j-1}$ gives the fundamental class
for $\CP^{j-1}$. (This calculus works readily for cone metrics
with differential forms that are suitably continuous.)
We conveniently obtain such a form as some constant multiple $\alpha$
of the K\"ahler form of the model geometry $\CP^{j-1}$ of the link.   
One way to determine $\alpha$ is to reduce to the case $j=2$ by
clustering the cone points into three groups which are collapsed along
a codimension 2 stratum limiting to $S$.  In the case $j=2$,
the complex link is $S^2$ with cone points of curvature
$\kappa_1+\kappa_2$, $\kappa_2+\kappa_3$ and $\kappa_3+\kappa_1$.
This uses up $2 \sum \kappa_i$ out of the total curvature
$4\pi$ of $S^2$, so the area of a constant curvature metric
is reduced by a factor of $\gamma(S)$.

Part (c) follows from (a) and (b), since the real link fraction
is the product of the complex link fraction with $\gamma(S)/2\pi$.

The case with symmetry follows by dividing the asymmetric configuration
space by the symmetry.
\end{proof}

\section {Orbifolds} 
\label{section: orbifolds}

An {\it orbifold} is a space locally modelled on $\reals^n$
modulo finite groups; the groups vary from point to point.  For
an exposition of the basic theory of orbifolds, see \cite{Thurston:GT3M}.
Our orbifolds will be $(X, G)$--orbifolds, locally
modelled on a homogeneous space $X$ with a Lie Group $G$ of isometries.
It is easily seen by induction on dimension
that an orientable $(X, G)$--orbifold has an induced metric which makes it into
a cone-manifold.  (Use the naturality of the exponential map.)

Here is a basic fact about the relation between
cone-manifolds and orbifolds, which essentially is a rephrasing
of Poincar\'e's theorem  on fundamental domains:

\begin{theorem}[Codimension 2 conditions suffice]
\label{codimension 2 conditions suffice}
Let $C$ be an $(X, G)$--cone-manifold.  Then $C$ is a ``weakening''
of the structure of an orbifold if and only
if all the codimension $2$ strata of $C$ have cone angles which
that are integral divisors of $2 \pi$.
\end{theorem}

\begin{proof}
An orientation-preserving group of isometries whose fixed point set
has codimension $2$ is a subgroup of $SO(2)$, and the only possibilities
are $\integers/n$.  The cone angle along such a stratum in an orbifold
is therefore an integral divisor of $2 \pi$, and the condition is
necessary.

The converse can be proved by induction on the codimension of the
singular strata of $C$.  Clearly, it works for strata of codimension $2$.
Suppose that we have proven that $C$ has an orbifold structure in
the neighborhood of all strata up through codimension $k$. 
Let $S$ be a singular stratum of codimension $k+1$, and consider
the neighborhood $U$ of a point $x \in S$.
This neighborhood can be taken to have the form of a bundle over a neighborhood
of $x$ in $S$, with fiber the cone on a $k$--dimensional cone-manifold $N$,
the normal sphere to $S$.  The normal sphere $N$ is modelled on $(S^k, G)$,
where $G \subset SO(n)$.  By induction, $N$ is an orbifold; its universal
cover must be $S^k$, since for $k \ge 2$ the sphere is simply-connected.
Therefore the cone on $N$ is the quotient of $B^{k+1}$ by the group of
covering transformations of $S^k$ over $N$, and therefore $U$ is also
the quotient space of a neighborhood in $X$ by the same group.  Thus,
$C$ is an orbifold.
\end{proof}

To illustrate, let's look at some of the local orbifold
structures that arise in multi-way collisions. When
$k$ cone points of equal curvature $2\pi \alpha$ collide,
the order of the local group $\Gamma(S)$
for a stratum $S$ is the reciprocal of the real volume fraction, so from
\ref{proposition: multiple cone angles}, setting $\alpha = \gamma(S)/2\pi$
 we have
\[ \#\Gamma(S) = \frac{k!}{ (1- k \alpha)^{k-1}}  \quad
\left [ \frac{1}{1/2 - \alpha} \in \integers  \quad\& \quad
0 < \alpha < 1/k \right ] \]
The only three cases satisfying the condition for three colliding equal
angle cone points are when $\alpha$ is $1/6$, $1/4$ and $3/10$.
The complex links in these
three cases are the quotient orbifolds of the sphere by the
oriented symmetries of one of the regular polyhedra: $(2,3,3)$,
$(2,3,4)$ or $(2,3,5)$. The real link is $S^3$ with a cone
axis along the trefoil knot of order 3, 4 or 5.  The formula
gives orders for these groups of 24, 96 and 600.  (This  can be 
quickly confirmed by an automated check
using the 3--dimensional topology program \texttt{Snappea}
to obtain presentations for the orbifold fundamental groups
and feeding then to a group theory program such as \texttt{Magnus}.)

An interesting example of a collision of cone points of unequal
curvatures is $(19\pi/30, 11\pi/30, 29\pi/30)$. The real
link is an orbifold with $(2,3,5)$ cone axes on the 3--component Hopf link.
In this case, $\alpha = 1/60$ and $\#(\Gamma(S)) = (60)^2 = 3600$.

The biggest possible multiple collision is when 5 points of curvature $\pi/3$
collide. The local group for this collision has order $6^4 5!= 155,520$.

Infinitely many of the modular spaces for cone-metrics with 4 cone points
are orbifolds of complex dimension 1, but for higher dimensional modular
spaces, only 94 are orbifolds.  These are tabulated in the appendix.

\section {Proof of main theorem} 

\proof[Proof of Theorem \ref{cone metrics form cone-manifold}]
Most of this theorem follows formally from
Theorem \ref{Completion is cone-manifold},
Proposition \ref{cone angles around collisions}, and
Theorem \ref{codimension 2 conditions suffice}.
What still remains is a discussion
of the volume of $\bar C(\alpha_1, \dots,
\alpha_k)$.

The only case in which
$X = \bar C(\alpha_1, \dots, \alpha_k)$ is not compact is where there are
cone-manifolds $x \in X$ whose diameters tend to infinity.  In such a case,
if we normalize so that the area of $x$ is $1$,
there must be subsets of $x$ with large diameter and small
area, free from cone points.  This implies that $x$ has subsets which
are approximately isometric to a thin Euclidean cylinder. If $\gamma \subset x$
is a short curve going around such an approximate cylinder, then the angle
of rotation for $\gamma$ must be a sum of a subset of the $\{\alpha_i\}$.
There are only a finite number of possibilities, so if the diameter
is large enough, a neighborhood of $\gamma$ of large diameter is
actually a cylinder.   Once $\gamma$ is determined,
the shapes of the two pieces of $x$ cut by $\gamma$ can be specified
independently, and a scale factor $\length(\gamma)^2/\area$
(less than some constant $\epsilon$)
together with an angle of rotation
can also be specified independently.

It will follow that the ends of $x$ are
in 1--1 correspondence with partitions $Q$ of the set of curvatures into
two subsets each summing to $2 \pi$, if we verify two points: 
\begin{itemize}
\item[\bf(i)]  for any such partition $Q$, there exists an $x \in X$ with a
geodesic $\gamma$ separating the cone points according to $Q$, and
\item[\bf(ii)] the subspace $X_{\gamma, \epsilon}$
consisting of cone-manifolds in
$X$ with a geodesic $\gamma$ of length $\epsilon$
which separates the cone points according to $Q$ is connected.
\end{itemize}

Actually, the proof does not logically depend on either point, and it
is a slight digression to prove them, but it seems worth doing anyway.

An easy demonstration of (i) is to construct a polygon with
angles $\pi - \alpha_i/2$.  It is easy to find a very thin
polygon realizing $Q$.  Doubling such a polygon gives a suitable
cone-manifold $x$.

We will describe an explicit construction for (ii).
Let us begin with the special case of $c \in X_{\gamma, \epsilon}$
which are obtained by doubling a convex Euclidean polygon whose angles
are half the cone angles for $X$.
It is easy to connect
any two convex polygons with the same sequence of angles
by a family of polygons having the same angles.
If we allow degenerate cases as well, where two angles coincide,
the order is irrelevant.  Therefore, this special subspace of
$X_{\gamma, \epsilon}$
is connected.  

Therefore, it suffices to connect any $c \in X_{\gamma, \epsilon}$
to something obtained by doubling a convex polygon.  Construct a
maximal cylindrical neighborhood $N_1$ of $\gamma$ with geodesic boundary.
There is at least one cone point on each boundary component of $N_1$.
Let $\beta$ be one of the boundary components, and $x_1 \in \beta$
a cone point, with curvature $\alpha$.  If $c$ is cut along $\beta$,
the portion on the other side of $\beta$ from $N_1$ has boundary
consisting of a geodesic with a convex angle of $\pi - \alpha$ at
$x_1$, and possibly additional angles if it contains other cone points.
There is a circular arc $\beta '$ through $x_1$,
contained in $N$, which appears to have a convex angle of $\pi - \alpha$
from within $N$, but appears to be smooth at $x_1$ when viewed from the outside.
Let $U_1$ be the ``outside'' component obtained by cutting along $\beta'$.
Its boundary is now locally isometric to a circle, and a neighborhood,
like on a cone, is foliated by parallel circles.

Deform $c$, by shrinking the ``interesting part'' of
$U_1$ relative to the rest of $c$, so that the next cone point in $U_1$ is not
close to $\beta'$.  Let $N_2$ be a maximal neighborhood of $\beta'$ which
is foliated by parallel circles, and let $x_2$ be a point on its boundary.
Adjust by a rotation of $U_1$ until the geodesic through
$x_1$ perpendicular to the foliation by circles hits at $x_2$.
Draw a circular arc through $x_2$, within $U_1$, which appears smooth
from the outside neighborhood $U_2$.

This process can be continued, in the same manner, until the last neighborhood
$U_k$ is a cone.  The geodesic through $x_{k-1}$ automatically hits the
cone point.  Now do the same process on the other side of $N_1$, first
adjusting by a rotation
so a geodesic through $x_1$ perpendicular to the foliation of
$N_1$ by parallel circles hits at a cone point.

After this sequence of deformations, we have
a cone-manifold with a geodesic Hamiltonian path through
all the cone points, such that at cone points internal to it the two
outgoing geodesics bisect the cone angle.  The path can be completed
to a curve by one additional geodesic (this is easy to see if you
draw the figure in the plane obtained by cutting along the path; it is
made of two convex arcs, and has bilateral symmetry).  

Note that a similar process works for a general cone-manifold: we do not
really need $\gamma$ for this construction, we can begin
at any cone point, and work outward from it.

We call the ends of $X$ {\it cusps}, in accordance with terminology
for manifolds and orbifolds.  To justify this word, note that each
cusp is foliated by complex geodesics with respect to the Hermitian
metric, obtained by rotating the two ends of $c$ with respect to each
other and by scaling.  The complex geodesics are locally isometric
to the hyperbolic plane.
The pure scaling, which may be thought of as 
inserting extra lengths of cylinder between the two ends,
generates a real geodesic.  These real
geodesics converge, as the shrinking increases.  The convergence
is exponential, so the total volume of each cusp is finite.
\endproof

\section {The icosahedron and other polyhedra}

Let $A$ be the subgroup of isometries of $\complexes$ which
take $\eisen$ into itself.
We may think of the classes of triangulations $P(n; k_1, \dots, k_m)$
as the space of $(\euclidean^2, A)$--cone-manifolds of area $n$ (measured
in double triangles) and cone angles $k_i \pi/3$. They consist of elements of
$C(k_1 \pi/3, \dots, k_n \pi/3)$ equipped with a reduction of the
$(\euclidean^2, \isom(\euclidean^2))$ structure to $(\euclidean^2, A)$.
In more concrete terms, a triangulation is given by a cocycle whose
coefficients are elements of $\eisen$.

Euclidean cone-manifolds sometimes admit several inequivalent
reductions to $(\euclidean^2, A)$---in other words, there are some
cone-manifolds that can be subdivided into unit equilateral
triangles in more than one way. 
In complex Lorentz space
$\complexes^{(m-3,1)}$,
the set of cocyles with a certain total area form a sheet of
a hyperboloid.  The hyperboloid fibers over complex hyperbolic space,
with fiber a circle (corresponding to multiplication of the cocycle by
elements of the unit circle).  The set of triangulations are lattice
points in $\complexes^{(m-3,1)}$, and the value of the Hermitian form
counts the number of triangles---multiple unit equilateral triangulations
of a Euclidean manifold correspond to fibers that intersect more than
lattice points. (All lattice points come in groups of 6 whose ratios are units
in the ring $\eisen$.)

\begin{figure}[htbp]
\begin{minipage}{.035\textwidth}\hbox to \textwidth {\hfill}
\end{minipage}
\begin{minipage}{.25\textwidth}
\includegraphics[width=\textwidth]{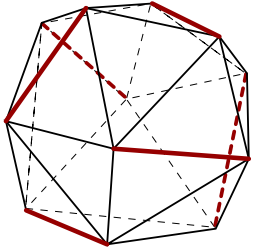}
\end{minipage}
\begin{minipage}{.045\textwidth}\hbox to \textwidth {\hfill}
\end{minipage}
\begin{minipage}{.65\textwidth}
\caption{If an icosahedron is slit along 6 disjoint arcs joining
its vertices in pairs, conical caps can be inserted to turn it into
an octahedron.
\label{figure: slit icosahedron}
}
\end{minipage}
\end{figure}

The ``biggest'' of the classes of triangulations is
\[
P(n; 1,1, \dots, 1) \subset J = C(\pi/3,\pi/3, \dots, \pi/3),
\]
the one which contains the icosahedraon.  The
``completion'' $\bar P(n; 1,\dots, 1) \subset \bar J$
which includes
degenerate cases contains all the other classes of triangulations.

By theorem 1.2, $\bar J$ is a complex hyperbolic orbifold
of dimension $9$.  The cone angles around the complex
codimension $1$ singular strata are $2 \pi/3$.

There is a concrete construction to describe an arbitrary element
of $\bar J$ or of $\bar P( 1, \dots, 1)$, as follows.
Suppose first that $x \in J$ is an arbitrary Euclidean
cone-metric on the sphere with all cone points having curvature $\pi/3$.
Choose a collection of $6$ disjoint geodesic arcs with endpoints on
the cone points.  Slit along each of these arcs.

Locally near the endpoints of the arcs, the developing map maps
the slit surface to the complement of a $60 \degrees$ angle.
A neighborhood of the slit develops to a region outside an equilateral
triangle in the plane; when you go once around the slit, the developing
image goes $2/3$ of the way around the triangle.

For each slit, take $2/3$ of an equilateral triangle with side
equal to the length of the slit, fold it together to form a cone
point in the center of the original triangle with curvature $2\pi/3$
and glue it into the slit.  The result is a cone-manifold $f(x)$
like the octahedron, in $C(2\pi/3, 2\pi/3, 2\pi/3, 2\pi/3, 2\pi/3, 2\pi/3)$.

As in section 2, we can analyze the shape of $f(x)$ by joining its cone
points in pairs by disjoint geodesic segments, slitting open,
and extending to give an element of $C(4\pi/3, 4\pi/3, 4\pi/3)$ (which is
a single point).

If $x \in \bar J - J$, the analysis still works: treat the
cone points as cone points with multiplicity, and use zero-length
slits as much as possible at cone points with curvature greater
than $\pi/3$.  At the first step, the slits of positive length pair
the cone points with curvature an odd multiple of $\pi/3$.  When the
slits are filled in, the curvature at each of the endpoints
is decreased by $\pi/3$, and the resulting cone-manifold has all curvatures
an even multiple of $\pi/3$.  For the second step, note that no cone
point can have curvature $6\pi/3$ or bigger. In this case, the slits of
positive length join cone points with curvature $2 \pi/3$.

An arbitrary $x \in \bar J$ can be reconstructed by reversing this
procedure.

\begin{figure}[htbp]
\centering
\begin{minipage}{.4\textwidth}
\centering
\includegraphics[width=.7\textwidth]{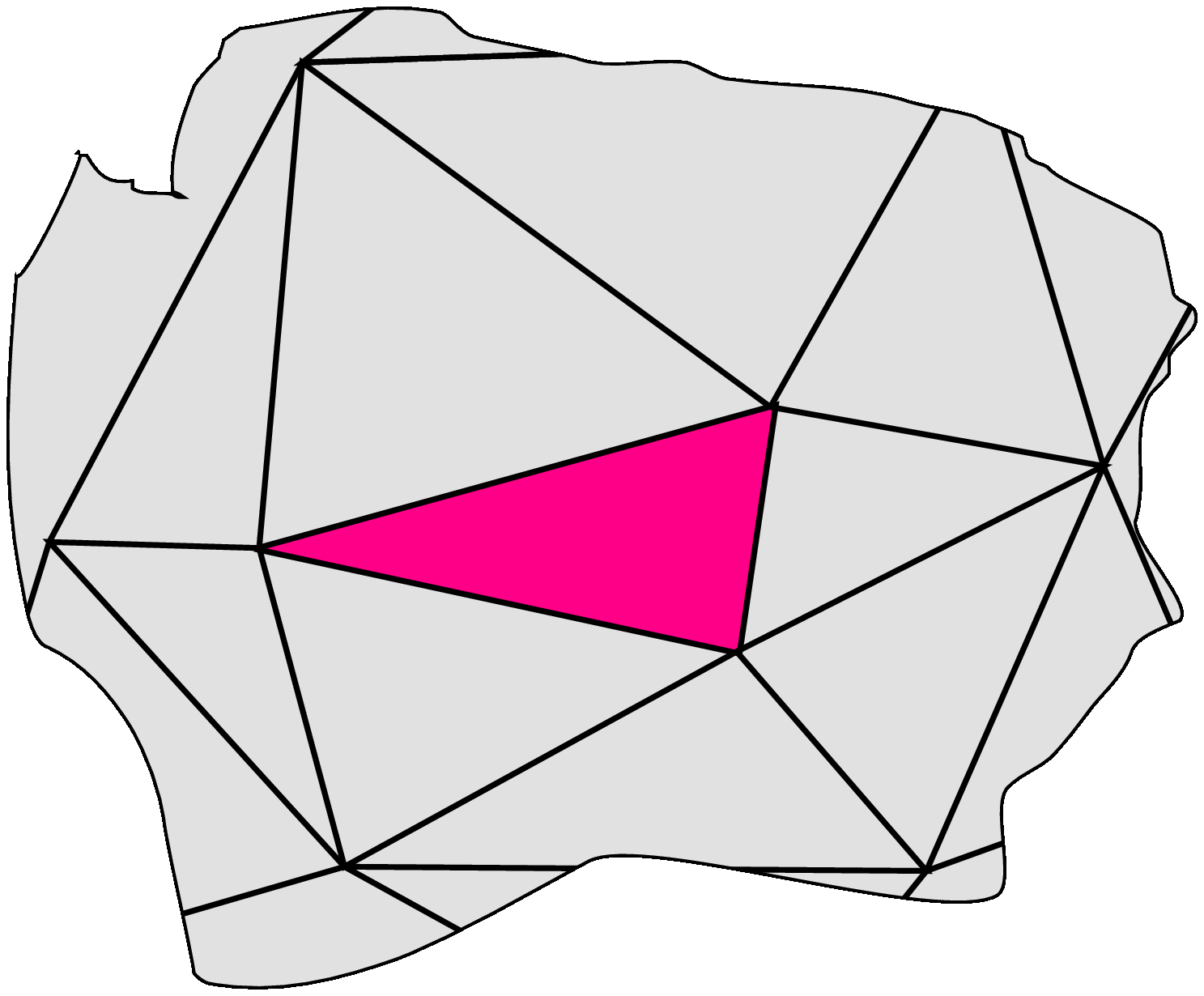}
\end{minipage}
\begin{minipage}{.4\textwidth}
\centering
\includegraphics[width=.6\textwidth]{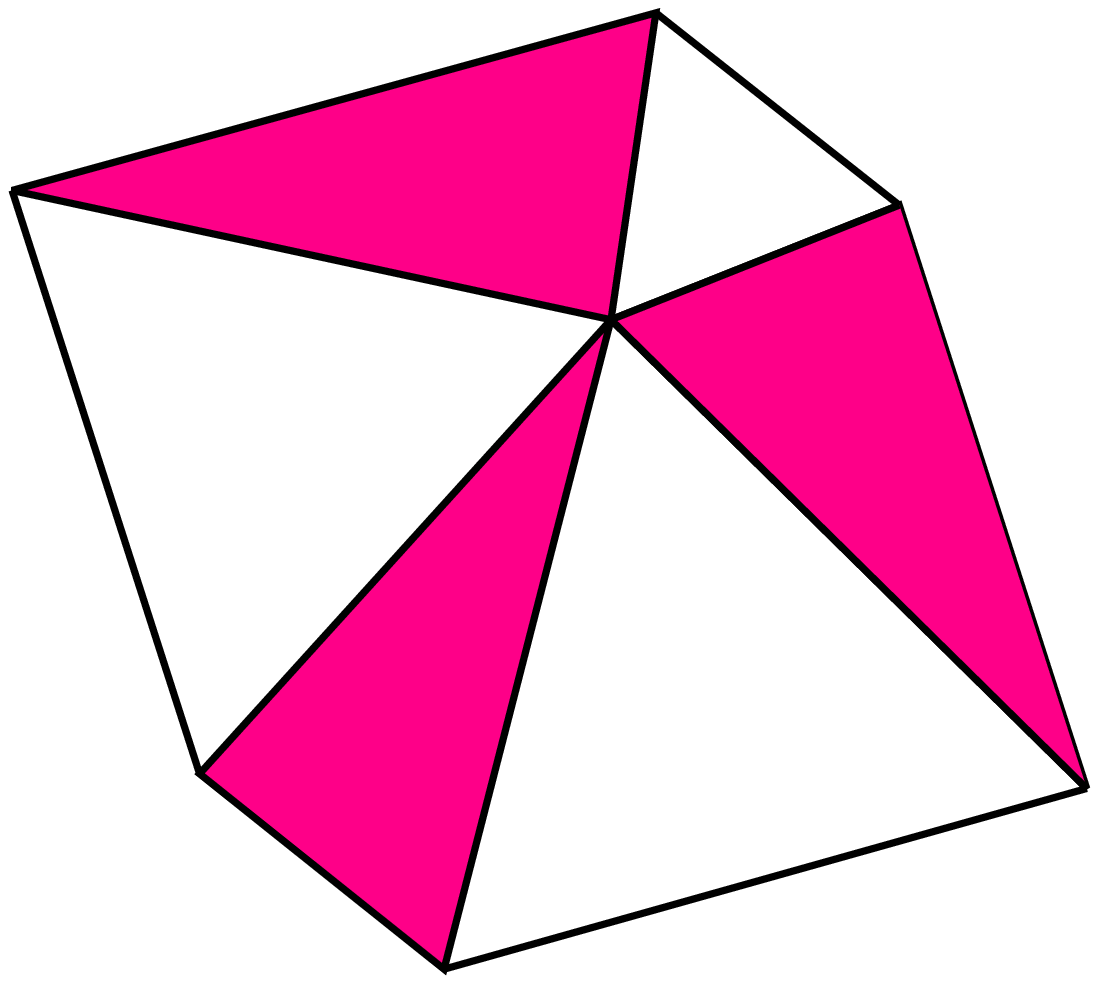}
\end{minipage}
\caption[Unrolling an icosahedron using Napoleon hexagons.]
{A geodesic triangle whose vertices are on 
cone points of curvature $\pi/3$ has a deleted neighborhood that
develops to the deleted neighborhood of a Napoleon hexagon,
formed from three copies of the triangle and three equilateral
triangles. This process, applied to the 12 vertices of
an icosahedron-like cone-manifold grouped into 3's, recursively
reduces it to a tetrahedron-like cone-manifold}
\end{figure}

There are many alternative coordinate systems for $J$.  For example,
another construction is to group the cone points in $3$'s, by constructing
$4$ disjoint geodesic triangles with vertices at cone points.  If
these triangles are cut out, then the developing image of what is
left is discrete; it comes from a $2,2,2,2$ group acting in
the Euclidean plane.  The developing image is the complement of
a certain union of hexagons about the lattice of elliptic points.
The hexagons are not arbitrary, however---the hexagons $H(T)$ that arise
are hexagons that come from Napoleon's theorem, constructed
as follows:
Suppose $T$ has sides $a$, $b$, and $c$, in counterclockwise order.
We will construct $6$ triangles around the vertex $v$ of $T$ between
$a$ and $b$.
First construct an equilateral triangle on side $a$.  Construct
another triangle $T_1$ congruent to $T$ on the free side of the
equilateral triangle which is incident to $v$.  Side $c$ of $T_1$
also touches $v$; on this, construct a second equilateral triangle.
Continue alternating copies of $T$ and equilateral triangles until
it closes, yielding $H(T)$.  

Note that $H(T)$ has sides $a, b, c, a, b, c$ in counterclockwise order.
The complement of $H(T)$ modulo a rotation of $180 \degrees$ has boundary
which matches the boundary of $T$; when it is glued in, three cone
points of curvature $\pi/3$ are obtained at the vertices of $T$.

A general $x \in \bar J$ can be obtained by choosing first a $2,2,2,2$
group, and then choosing four hexagons $H(T_i)$ centered about
the four classes of vertices.  Form the quotient of the complement
of the hexagons by the group, and glue in the triangles $T_i$.
If the hexagons are disjoint and nondegenerate, $x \in J$.

From this concrete point of view, what is amazing is that these
coordinate systems have a global meaning, since $\bar J$
is an orbifold:  even if one chooses
a collection of hexagons $H(T_i)$ which overlap, they determine
a unique Euclidean cone-manifold, provided the net area (computed
formally) is positive.

Using these constructions, it is not hard to show that 
$P(n; 1, \dots, 1)$ contains $1$ or more elements for all values of $n$
starting with $10$, with $11$ as the sole exception.  If there were
an element $T$ of $P(11, 1, \dots, 1)$, it would have $13$ vertices and
$22$ triangles.  One could then construct a spherical cone-manifold
by using equilateral spherical triangles with angles $2\pi/5$.  This
cone-manifold would have only one cone point --- which is manifestly
impossible, since the holonomy for a curve going
around the cone point is a rotation
of order $5$, but at the same time the holonomy is trivial since the curve
is the boundary of a disk having a spherical structure.

From the picture in $C^{(1,9)}$, it follows
that the number of non-negatively curved
triangulations having up to $2n$ triangles is roughly proportional to
the volume of the intersection of some cone with the ball of
radius $\sqrt n$ in this indefinite metric. The cone in question is
neither compact nor convex, but since it comes from a fundamental domain for
the group action, its intersection with the ball of norm less than any
constant has finite 10--real-dimensional volume.  Therefore, the
number of triangulations with up to $2n$ triangles is $O(n^{10})$.

\section{An explicit construction and fundamental domain}
\label{section: construction and fundamental domain}

Another method for constructing,
manipulating and analyzing non-negatively curved cone structures 
goes as follows:

\begin{description} \label{Construction: polygon to 
polyhedron}
	
		\item[Given] $k+1$ real numbers $\alpha_{0}, \alpha_{1}, \dots, 
	\alpha_{k} \ge 0$ whose sum is $4\pi$.
	\item[To Construct] Euclidean cone-metrics with the 
	$\alpha_{i}$ as 
	curvatures.
	
	\item[Choose] a $k$--gon $P$ in the plane, with edges 
	$e_{1}\dots, e_{k}$.
	
\item[Construct] ($i=1, \dots k$): An isosceles
	triangles $T_{i}$ with base on $e_{i}$, apex $v_{i}$,
	apex angle $\alpha_{i}$, pointing inward
	if $\alpha_{i} < \pi$, pointing outward if $\alpha_{i}>\pi$.
	
	\item[Condition A] the triangles $T_{i}$ are disjoint 
	from each other and disjoint from $P$ except along $e_{i}$.
	
	\item[Let] $Q$ be thefilled polygon obtained from $P$ by 
	replacing each $e_{i}$ by the other two sides $f_{i}$ and 
	$g_{i}$ of $T_{i}$.  
	
	\item[Glue] $f_{i}$ to $g_{i}$ to obtain a cone manifold.  
	The vertex $v_{i}$ becomes a cone point of curvature 
	$\alpha_{i}$.  The other $k$ vertices of $P$ all join to 
	form a cone point of cone angle $\alpha_{0}$.

\end{description}

\begin{figure}[hbtp]
\cl{\includegraphics[width=.4\textwidth]{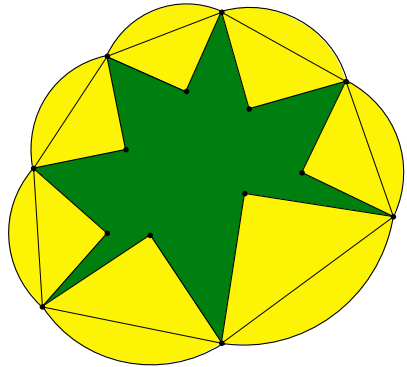}}
	\caption{A cube-like cone-metric (8 cone-angles of 
	curvature $\pi/2$) can be 
	constructed by removing isosceles 
right triangles from the sides of a heptagon and
gluing the resulting pairs of equal sides. The seven 
sharp angles all come together to form  the eigth cone point.
This illustration (along with the others in this section) was constructed 
with the program {\em Geometer's 
Sketchpad\/}, where the shape can be varied while preserving 
the correct geometric relations.
\label{figure: cuboid island}
}
\end{figure}

As examples, see figure \ref{figure: cuboid island} for 
a cube-like cone-manifold, or figure \ref{figure: 23 vertices} for a 
triangulation of 
$S^{2}$ with 23 vertices and 42 triangles constructed from an 
icosahedral-like cone-manifold.

Here is the inverse construction.  Given a cone-metric with $n$ cone points on 
$S^{2}$:
\begin{description}
	\item[Choose] one of the 
cone points $v_{0}$.
\item[Find] for each other cone point $v_{i}$ a shortest path
$a_{i}$ from $v_{0}$ to $v_{i}$.
The $a_{i}$ are necessarily simple and 
disjoint, except at $v_{0}$.
\item[Cut]  along all these 
paths, to obtain a disk equipped with a Euclidean metric 
whose boundary is composed of $2 (n-1)$ straight segments, 
each paired to an adjacent segment of the same length and 
forming an angle equal to the corresponding cone angle. (See 
figure \ref{figure: icosahedral star}.)
\end{description}

We will show that if $P$ is a cone-metric on the sphere with positive 
curvature at each vertex, and if $S(P)$ ($S$ because it 
resembles a star) is the metric on 
$D^{2}$ obtained by cutting $P$ open as above, then $S(P)$ 
can be flattened out into the plane, that is, it is 
isometric to the metric of a filled simple polygon.\break

\cl{\includegraphics[width=4.2in]{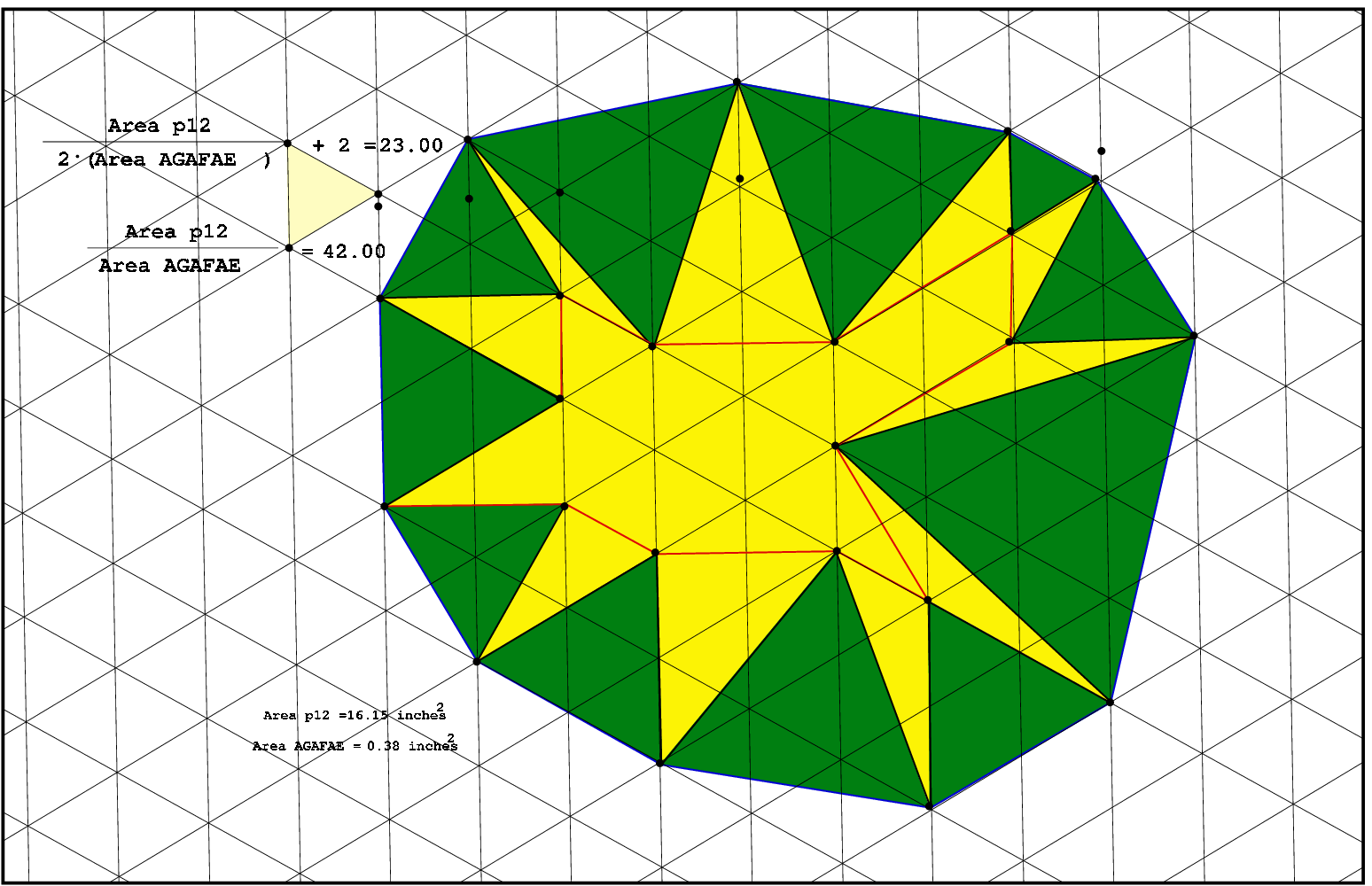}}\vglue -5mm
\begin{figure}[hbtp]	
\caption{This is a diagram for a triangulation of $S^{2}$ 
	with 42 triangles, having 12 vertices of order 5 and 11 
	vertices of order 6. Eleven--gons whose vertices lie in an 
	ideal of index 3 (generated by $1 \pm \omega$) in the Eisenstein 
	integers determine non-negatively-curved triangulations of 
	$S^{2}$. Each valley between star-tips is folded together 
	to form the triangulation; the star-tips come together at 
	the base vertex.
	If the inner vertices of the 11--pointed star are 
	closer to the two adjacent star-tips than to any other 
	star-tips, this is the canonical 11--gon for the 
	triangulation based at the given vertex.
	\label{figure: 23 vertices}
	}
\end{figure}

Actually, we will enlarge $S(P)$ to a surface $F(P)$ 
(resembling a flower) by adjoining
sectors of circles of  
with apex at 
each vertex $v_{i} \, (i > 0)$ of $S(P)$ and angle equal to 
the curvature at $v_{i}$ in $P$, so that the resulting 
surface is locally Euclidean everywhere in its interior (as 
in figure \ref{figure: cuboid island}).

The minimum distance within $S(P)$ of any point in $S(P)$ from one 
star-points that assemble at $v_{0}$ is equal to the 
distance of its image in $P$ from $v_{0}$. Let $Q \subset 
S(P) \subset F(P)$ be the set of points whose minimum 
distance to $\partial F(P)$ is attained at 3 or 
more points on $\partial F(P)$. Then $Q$ includes 
$\set{v_{1}, \dots, v_{n}}$, as well as the  vertices 
for the Voronoi diagram of the  star tips within $S(P)$. 
Let $R$ be the collection of open segments consisting of points 
whose minimum distance to $\partial F(P)$ is attained at 
two points of $\partial F(P)$; they are the edges of a tree, 
whose vertex set is $Q$. We
decompose $S(P)$ into dart-shaped quadrilaterals, consisting of 
union of the two minimum-length arcs from points in an 
edge $\alpha \in R$
 to $\partial F(P)$ (see figure \ref{figure: 
icosahedral star}). We'll call this dart $D(\alpha)$.  Let 
$\theta(e)$ be the angle of $D(e)$ at either of its two 
wingtips (vertices that are not vertics of $e$).  
Note that 
\[
\sum_{e \in R} \theta(e) = 1/2 \left ( 2 \pi - 
\kappa(v_{0}) \right ),
\]

\cl{\includegraphics[width=.6\textwidth]{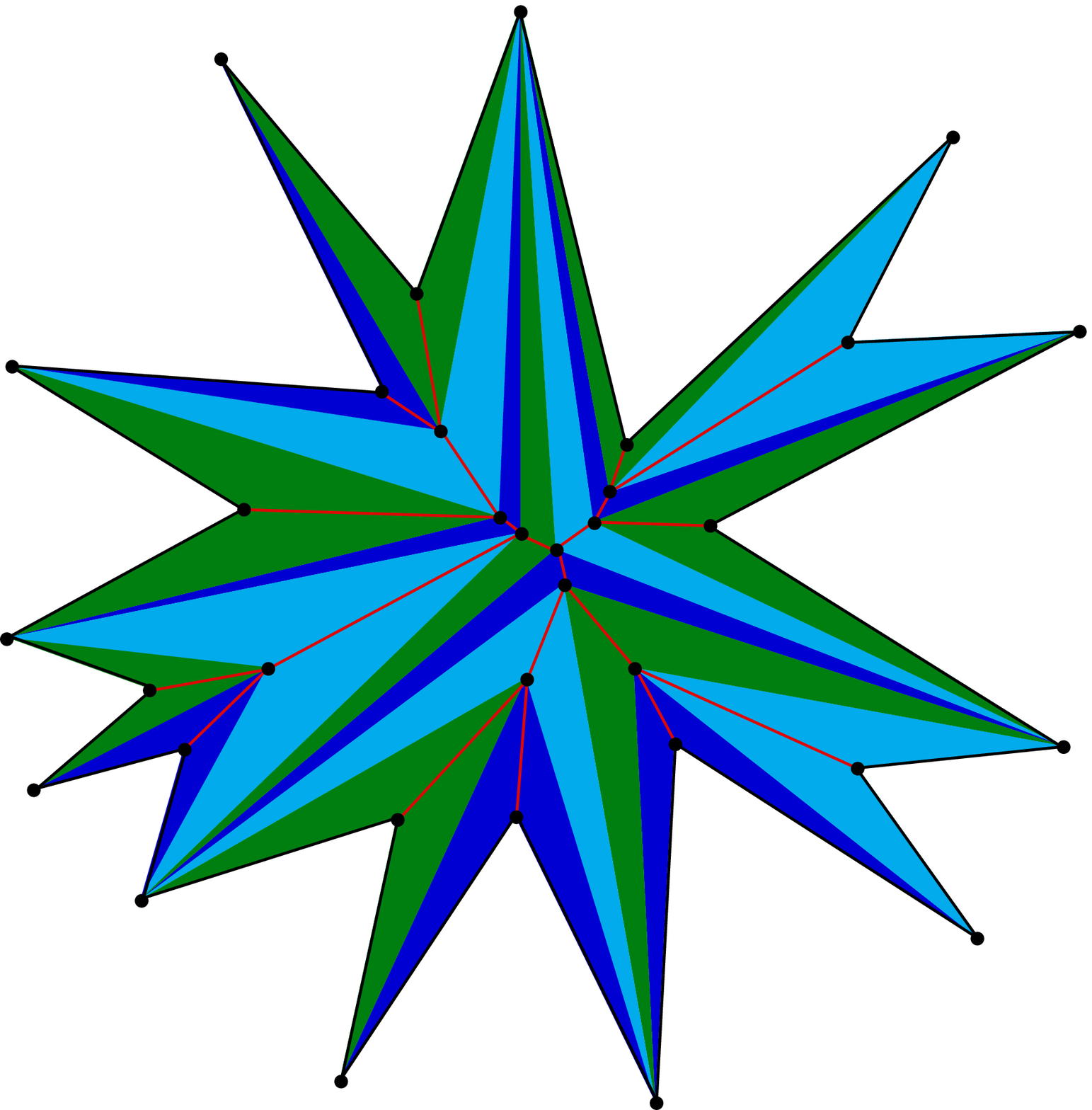}}\vglue -6mm
\begin{figure}[htbp]
	\caption{An irregular icosahedron sliced and flattened. A
regular point $\bigstar$ was chosen on a 
	Euclidean cone-metric for $S^2$ having 12 cone points of 
curvature $\pi/3$. The surface has been cut along shortest geodesics
from 
	$\bigstar$ to each of the cone points, and flattened 
	into the plane to form a 24--gon resembling a star. The polygon has been
	subdivided into 45 dart-shaped
	 quadrilaterals. Each quadrilateral is obtained from an 
	 edge of the cut locus of the original icosahedron ($=$ 
	 the Voronoi diagram for the $\bigstar$--tips, after cutting)  by 
	 suspending to its two closest $\bigstar$--tips. \label{figure: 
	 icosahedral star}
	 }
\end{figure}

that is, half the cone angle at the base vertex.

For any vertex $q \in Q$, let $D(q) \subset F(P)$ be the 
maximal disk in $F(P)$ centered about 
$q_{i}$. If $q_{1}$ and $q_{2}$ are the endpoints of 
$e \in R$, then the angle between the bounding circles 
of $D(q_{1})$ and $D(q_{2})$ is $\theta(e)$.   

\begin{proposition}
	$F(P)$ has an isometric embedding in the plane $\euclidean^{2}$.
\end{proposition}
\begin{proof}
The developing map $f\co  F(P) \to \euclidean^{2}$ into the plane is  an isometric 
immersion. To show that it is an embedding, it will suffice to 
establish that $f$ restricted to the boundary $\partial F(P)$ 
is an embedding.

The boundary $\partial F(P)$ is composed of inward-curving 
circle arcs that meet at outward-bending angles.  For each 
edge $e \in R$, there is a pair of these angles, where 
$\partial F(P)$ 
turns 
by an angle $\theta(e)$.  For any two points $x, y \in 
\partial F(P)$, there is at least one path along $\partial 
F(P)$ where these bending angles sum to no more than $\pi$.

\cl{\includegraphics[width=.6\textwidth]{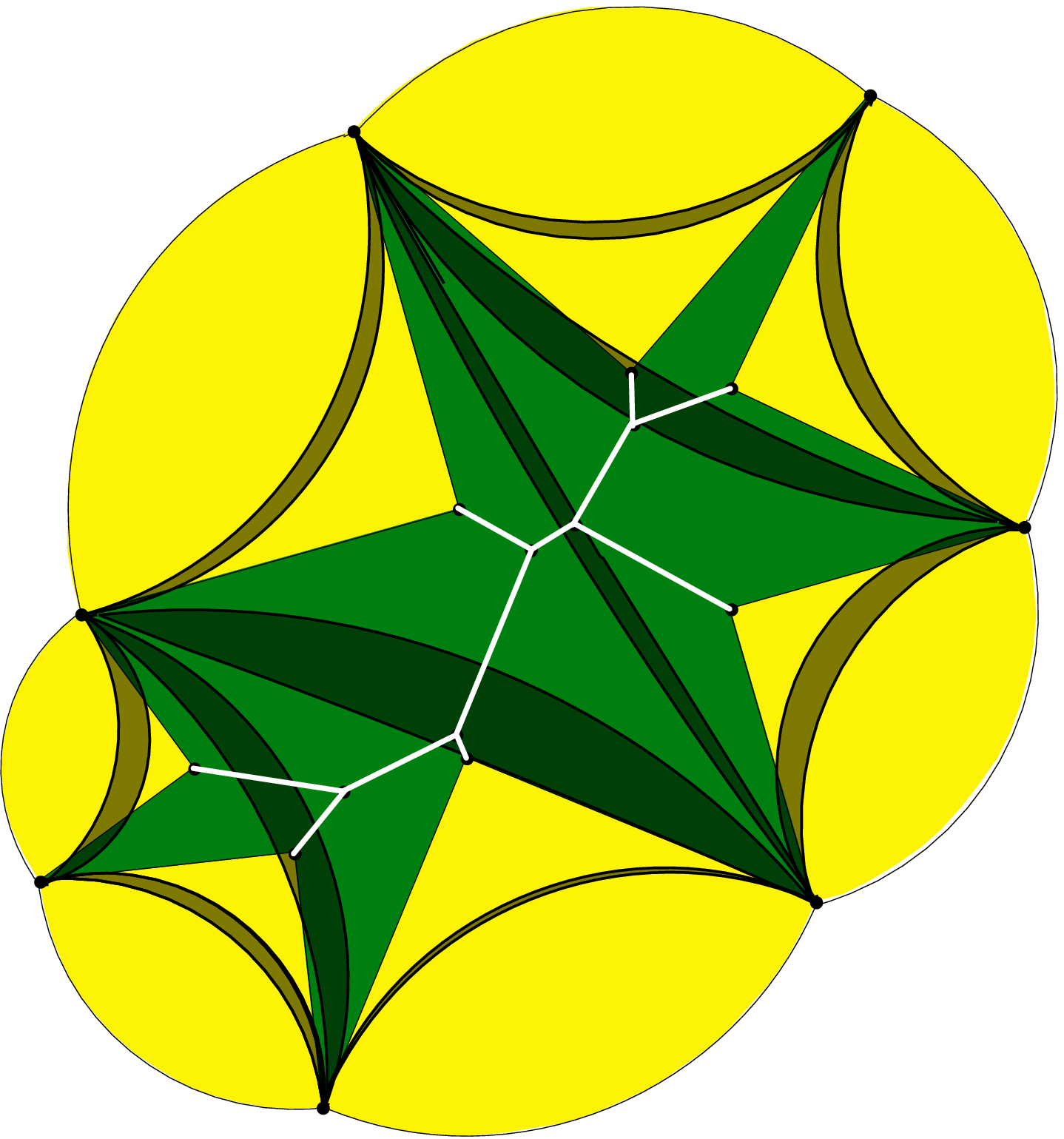}}
\vglue -6mm
\begin{figure}[htbp]
	\caption{A hyperbolic view associated with the 
	cut-open polyhedron.  From the point of view of 
	3--dimensional hyperbolic geometry, if this figure is 
	interpreted as lying on the boundary of upper half-space, 
	the convex hull of its complement is the union of the 
	hemispherical bubbles which rest on it. The boundary of 
the convex hull (with geometry induced from the upper half-space metric 
$ds^{2} = 1/z^{2}(dx^{2}+dy^{2}+dz^{2})$) is isometric with the hyperbolic
	plane, bent into hyperbolic 3--space.  The sum of all 
	bending angles is one half the cone angle at the base 
	point (assembled from the tips of stars). Any immersion of 
	the hyperbolic plane which has total bending measure less 
	than $\pi$ is an embedding, so this plane is embedded.  
	There are immersed planes with total bending any number 
	greater than $\pi$ which are not embedded.
	\label{figure: cuboid poincare bending}
	}
	\end{figure}

For an immersion of a disk in $\euclidean^{2}$ to fail to be an 
embedding, any innermost arc whose endpoints are identified 
by the immersion must have total curvature at least $-\pi$ 
when orientations are chosen so that the total curvature of 
the entire boundary is $2 \pi$.  This is clearly impossible 
in our situation, so $f$ is actually an embedding.
\end{proof}

\begin{remark}
	This proposition can be rephrased in terms of
	3--dimensional hyperbolic geometry: any pleated immersion
	$f\co \hy^{2} \to \hy^{3}$ with positive bending measure whose 
	integral on any geodesic is no greater than $\pi$ is an 
	embedding. This is related to the inequality of Sullivan 
	analyzed and refined by Epstein and Marden in 
	\cite{MR89c:52014},
	and also to the global characterization of bending for 
	convex polyhedra of Rivin and Hodgson \cite{MR93j:52015}, 
	see Rivin \cite{MR96i:52008} for a related inequality for
	convex hyperbolic polyhedra.
\end{remark}

The $n$--dimensional \df{associahedron} is a polyhedron 
whose vertices are labelled by triangulations of an 
$(n+3)$--gon using only vertices of the polygon, and whose 
$k$--cells are labelled by subdivisions of the $(n+3)$--gon 
obtained by removing $k$ edges from a triangulation. They 
can be thought of as describing all ways to parenthesize or 
associate a product of $n+2$ symbols. The associahedron is a convex 
polyhedron in $n$--space that arises in a variety of 
mathematical context including the theory of 
loop spaces, Teichm\"uller theory and numerous 
combinatorial settings.  The numbers of triangulations are 
called Catalan numbers.

If cone angles at the $n$ points of $P$ are fixed,
the angles $\theta(e)$ that can occur for our dart quadrilaterals
can be described by mapping the $n-1$ star tips 
of $S(P)$ to the vertices of a regular $(n-1)$--gon, mapping 
each dart quadrilateral $D(e)$ to an edge or a chord of 
this polygon, and labelling each edge by the angle $\theta(e)$.
(In terms of hyperbolic geometry, this is an element of the 
measured
lamination space for the ideal polygonal orbifold $(\infty, \infty, \dots, 
\infty)$.)  This
determines a point in a polyhedron 
$\mathcal{F}n$
closely related to an 
associahedron, namely, the join of the boundary of the 
dual of the $(n-4)$--dimensional 
associahedron with the $(n-2)$--simplex. (When the measure on 
the boundary of the polygon is zero, we get a point on the 
boundary of the dual of the $(n-4)$--dimensional associahedron.  
Measures on the polygon itself with fixed total weight form an 
$(n-2)$--simplex.)

The set of all  possible functions  $\theta(e)$ (which we refer to as 
measures, after the usage in hyperbolic geometry and 
Teichm\"uller theory) can be described globally 
as a convex polyhedron using dual train track coordinates, 
as follows:  rotate a copy of the regular $(n-1)$--gon $1/(2n-2)$th 
of a revolution so it is out of phase 
with itself.  Choose any triangulation of this rotated 
polygon, using only its vertices. For each edge $f$ of this 
triangulation, let $m(f)$ be the sum of $\theta(e)$ where 
$e$ intersects $f$.  For any triangle with sides $f,g,h$, 
the
quantities $m(f)$, $m(g)$ and $m(h)$ satisfy the three triangle 
inequalities $m(f) + m(g) \ge m(h)$ etc.  These 
measures are subject to one linear constraint, namely, the 
sum of $m(f)$ where $f$ ranges over the edges of the 
rotated polygon adds to the cone angle at the base vertex 
$v_{0}$. 

For any set of numbers $\set{m(f)}$ satisfying 
the linear equation and linear inequalities, a measured 
lamination having total measure $\pi -\alpha_{0}/2$, where 
$\alpha_{0}$ is the cone angle at $v_{0}$,
can be reconstructed by a simple method familiar in the 
theory of measured foliations or normal curves on surfaces,
 by first solving for the picture in 
each triangle of the rotated polygon, then gluing the triangles together.
From this measured lamination and from the specification of 
cone angles (in order) at $v_{1}, \dots v_{n-1}$, a star 
polygon in the plane can be constructed recursively, using 
the principle that the shape of any dart quadrilateral 
$D(e)$ is determined from $\theta(e)$ together with either 
of its other two angles.  This star polygon is determined 
up to similarity.  When glued together 
it forms a cone manifold with specified cone-angles.

If all cone angles are equal, and if we are not 
distinguishing shapes that are the same up to permutation of 
the labels of cone points $v_{1}, \dots v_{n-1}$, then we must divide $\mathcal{F}$ 
by action of the group of order $n-1$ rotations. The 
faces of $\mathcal F$ correspond to measures $\theta$ where one of
the edges of the $(n-1)$--gon has measure 0.  Geometrically, 
this means that one of the cone points $v_{i}$, $i > 0$ has 
two or more shortest paths on $P$ to $v_{0}$. We could cut 
$P$ open along either of these shortest paths. In $S(P)$, 
this means one of the ``inside'' vertices of the star has 
three or more shortest paths to the tip vertices: two are 
sides of $S(P)$, and at least one is interior to $S(P)$. 
You can cut $S(P)$ along such an edge, and rotate one 
resulting chunk with respect to the other, to obtain a new 
shape $S'(P)$ with vertices in a permuted order.  

To also
insist on allowing change of base point requires 
a further much less direct equivalence relation.  
If the cone angles $\alpha_{1}, \dots, \alpha_{n-1}$ are 
not all the same, then to get all possible cone-metrics, we 
need one copy of $\mathcal{F}$ for each ordering of the cone 
angles up to cyclic permutation.

\section {Teichm\"uller space interpretation} 

Each element of $C(\alpha_1, \dots, \alpha_n)$ determines a point
in a certain finite sheeted covering of the modular orbifold for the
$n$--punctured sphere.  (The covering corresponds to the subgroup of
the modular group for the $n$--punctured sphere which preserves the
cone angles):  the map consists of forgetting the metric, and
remembering only the conformal structure.

By the uniformization theorem,
each of these metrics is equivalent to a metric obtained by deleting
$n$ points from the Riemann sphere $\Chat$. The resulting configuration
of $n$ points in $\Chat$ is unique up to M\"obius transformations.

\begin{proposition}
\label{cone metrics parametrize Teichmuller space}
The map from $C(\alpha_1, \dots, \alpha_n)$ is a homeomorphism.
\end{proposition}

\begin{proof}
In fact, there is an explicit formula for the inverse map,
going from a configuration of $n$ points on $\Chat$ together with
the curvatures at those points to a Euclidean cone-manifold with the
given conformal structure.  The formula is essentially
the same as the Schwarz--Christoffel
formula for uniformizing a rectilinear polygon. (See 
\cite{MR91h:53059} for an analysis of these and other
cone-manifold structures.)

The idea is to think of the construction of a Euclidean cone metric on $\Chat$
in terms of its developing map $h$.  Consider a collection
$\{y_i\}$ of
points in $\Chat$, with desired curvatures $\{\alpha_i\}$.  Let $P$ be the
punctured Riemann sphere $\Chat - \{y_i\}$.
The developing map $h$ is not uniquely
determined on $P$, and it is only defined on the universal cover $\tilde P$,
but any two choices differ by
a complex affine transformation.  Therefore, the {\it pre-Schwarzian}
of $h$, that is, $S =  h''/h'$, is uniquely determined by the metric, and
it is defined on $P$, not just on the universal cover of $P$.
The Euclidean metric can be easily reconstructed if we are given $S$, because
once we choose an initial value and derivative for the developing map $h$
at one point on $\tilde P$, we can integrate the differential equation
$h'' = S h'$ to determine it everywhere else.

How can we determine $S$?  Consider a cone, with curvature $\alpha$ at
the its apex.  If a cone is conformally mapped to $\complexes$
with its apex going to the origin, the developing
map is $z \gozto z^{1-\frac{\alpha}{2\pi}}$.  The pre-Schwarzian for this
map is $-\frac{\alpha}{2\pi} z^{-1}$.  It follows that the pre-Schwarzian
for the developing map of any Euclidean cone-metric with a cone point
having curvature $\alpha$ will have a pole at the cone point, with
residue $-\frac{\alpha}{2\pi} z^{-1}$.  Conversely, if the pre-Schwarzian
of some function $h$
has a pole of this type at any point in $\Chat$, then $h$ will locally
be the developing map for a Euclidean structure with a cone point of angle
$\alpha$.  (To see this, observe that the analytic continuation of $h$
around the pole differs by post-composition with an affine transformation.
Using this information, one can make a local conformal change of coordinates
in the domain so that $h$ has the form $z \gozto z^{1-\frac{\alpha}{2\pi}}$,
where $\alpha$ is not necessarily real.  From this, one sees that the
pre-Schwarzian has a pole at the singularity with residue $\alpha/{2\pi}$.)

We may as well assume that the the $\{y_i\}$ are in the finite part of
$\Chat$.  Define
$$S = \sum_i - \frac{\alpha}{2 \pi} (z-y^{-1}) . $$
Computation shows that in a coordinate patch $w = z^{-1}$
for a neighborhood of $\infty$, the pre-Schwarzian in terms of the variable
$w$ is holomorphic
if and only if $S$ behaves asymptotically like $-2 z^{-1}$.  This is
satisfied in our case, since the sum of the $\alpha_i$ is $4 \pi$.
The condition that $S$ is holomorphic on $P$, and that it has the
given behaviour at the cone points and at $\infty$, uniquely determines
$S$.

$S$ determines a complex affine structure on $P$.  Since the
fundamental group of $P$ is generated by loops going around the punctures,
and since the holonomy around these loops is isometric, the affine structure
is compatible with a Euclidean structure, well-defined up to scaling.
\end{proof}

Thus, we may think of $C(\alpha_1, \dots, \alpha_n)$ as a certain
geometric interpretation of modular space.
The completions $\bar C(\alpha_1, \dots, \alpha_n)$ have a topology
which depends on the comparisons
of sums of subsets of the $\alpha_i$ with $2 \pi$.  
It is almost never agrees with the standard compactification of the modular
space.  However, there are only a finite number of possible possibilities
for the topology --- it is curious that
we thus obtain several parameter families of complex hyperbolic structures on
the Teichm\"uller space, and several parameter families of complex hyperbolic 
cone-manifolds on the various $\bar C(\alpha_1, \dots, \alpha_n)$, with
varying cone angles.

Is there any similar phenomenon for the Teichm\"uller spaces of 
other surfaces, particularly closed surfaces?  The surface of genus $2$
has the same modular space as the six--punctured sphere, so for that
particular case, the construction carries over. I don't know how to
extend it to surfaces of higher genus.

\Addresses\recd

\small

\section*{Appendix: 94 orbifolds}

We give below a list of the examples of the spaces
$C(\alpha_1, \dots, \alpha_n)$ which are orbifolds, for $n \ge 5$.
When $n=3$ there is only one example for each feasible triple of cone angles,
and for $n = 4$ there are infinitely many examples.
In fact,
every triangle group in the hyperbolic plane can be interpreted as the modular
space for families of tetrahedra.  In general, the
$\alpha_i$ are
of the form $\frac{2\pi p}{q}$, for $p$ and $q$ integers.  For each
example, we list the least
denominator $q$ and the sequence of numerators $p_i$.  We also list
the degree of the number field containing the roots of unity
$\exp(\frac{2\pi p_i}{q})$ (that is, the number of integers less than
$q$ relatively prime to $q$).
We list also three
additional bits of information:
\begin{description}
\item[(arithmetic)]  Is the orbifold arithmetic (AR) or non-arithmetic (NR)?
\item[(pure)]  Is the completion of the covering of the modular  space whose
fundamental group is the pure braid group an orbifold (P),
or are some interchanges
of pairs of cone points needed to make the orbifold (I)?
\item[(compact)]  Is the orbifold compact (C) or non-compact (N)?
\end{description}

The question of arithmeticity hinges on the signatures of the Hermitian
forms obtained when we conjugate the curvatures at the cone points
(considered as roots of unity) by the Galois automorphisms.  If all the
other signatures are negative definite or positive definite, the group
is arithmetic; otherwise not.  The other two questions are more obvious.

These examples were enumerated by a routine computer program, which checks
all possibilities having a given least common denominator $q$.  The
enumeration was not rigorously verified (even though it should not be hard
to do so and search more `intelligently' at the same time)
but was a simple check of all denominators through $999$ in a few minutes
of computer time.
Mostow has rigorously enumerated examples by hand,
so this table can be regarded as just a check.\newpage

\font\table=cmr8

\vbox{\offinterlineskip\table
\def\marker{\vrule height 0.122truein depth 0.03truein width 0.7pt}
\def\topmarker{\vrule height 0.12truein depth 0.04truein width 0.7pt}
\def\tablerule{\noalign{\hrule height0.7pt}}
\let\\\cr
\def\s{\hglue 3.6pt}
\tabskip=0pt
\halign{\marker\hfil\s#\s\hfil&\marker\hfil\s#\s\hfil
&\marker\hfil\s#\s\hfil&\hfil\s#\s\hfil&\hfil\s#\s\hfil&\hfil\s#\s\hfil
&\hfil\s#\s\hfil&\hfil\s#\s\hfil&\hfil\s#\s\hfil&\hfil\s#\s\hfil
&\hfil\s#\s\hfil&\hfil\s#\s\hfil&\hfil\s#\s\hfil&\hfil\s#\s\hfil
&\marker\hfil\s#\s\hfil&\marker\hfil\s#\s\hfil
&\marker\hfil\s#\s\hfil&\marker\hfil\s#\s\hfil\marker\cr
\tablerule
\omit\topmarker&\omit\topmarker\s{\tiny Denominator}\s&\omit\topmarker\hfil
&\multispan 5
\hfil{\tiny Numerators}\hfil&&&&&&&\omit\topmarker\s{\tiny degree}\s
&\omit\topmarker\s{\tiny arithmetic?}\s
&\omit\topmarker\s{\tiny pure?}\s
&\omit\topmarker\s{\tiny compact?}\s\topmarker\cr
\tablerule
\nex &\quad 3  & 1 &  1 &  1 &  1 &  1 &  1 &&&&&&&\quad 2 &AR &P &N\\
\nex&\quad 3  &  2 &  1 &  1 &  1 &  1   &&&&&&&&\quad 2  &AR &P &N \\
\nex&\quad 4  &  1 &  1 &  1 &  1 &  1 &  1 &  1 &  1   &&&&&\quad 2  &AR &P &N \cr
\nex&\quad 4  &  2 &  1 &  1 &  1 &  1 &  1 &  1   &&&&&&\quad 2  &AR &P &N \\
\nex&\quad 4  &  3 &  1 &  1 &  1 &  1 &  1   &&&&&&&\quad 2  &AR &P &N \\
\nex&\quad 4  &  2 &  2 &  1 &  1 &  1 &  1   &&&&&&&\quad 2  &AR &P &N \\
\nex&\quad 4  &  3 &  2 &  1 &  1 &  1   &&&&&&&&\quad 2  &AR &P &N \\
\nex&\quad 4  &  2 &  2 &  2 &  1 &  1   &&&&&&&&\quad 2  &AR &P &N \\
\nex&\quad 5  &  2 &  2 &  2 &  2 &  2   &&&&&&&&\quad 4  &AR &P &C \\
\nex&\quad 6  &  1 &  1 &  1 &  1 &  1 &  1 &  1 &  1 &  1 &  1 &  1 &1&\quad 2  &AR &I &N \\
\nex&\quad 6  &  2 &  1 &  1 &  1 &  1 &  1 &  1 &  1 &  1 &  1 &  1   &&\quad 2  &AR &I &N \\
\nex&\quad 6  &  3 &  1 &  1 &  1 &  1 &  1 &  1 &  1 &  1 &  1   &&&\quad 2  &AR &I &N \\
\nex&\quad 6  &  2 &  2 &  1 &  1 &  1 &  1 &  1 &  1 &  1 &  1   &&&\quad 2  &AR &I &N \\
\nex&\quad 6  &  4 &  1 &  1 &  1 &  1 &  1 &  1 &  1 &  1   &&&&\quad 2  &AR &I &N \\
\nex&\quad 6  &  3 &  2 &  1 &  1 &  1 &  1 &  1 &  1 &  1   &&&&\quad 2  &AR &I &N \\
\nex&\quad 6  &  5 &  1 &  1 &  1 &  1 &  1 &  1 &  1   &&&&&\quad 2  &AR &I &N \\
\nex&\quad 6  &  2 &  2 &  2 &  1 &  1 &  1 &  1 &  1 &  1   &&&&\quad 2  &AR &I &N \\
\nex&\quad 6  &  4 &  2 &  1 &  1 &  1 &  1 &  1 &  1   &&&&&\quad 2  &AR &I &N \\
\nex&\quad 6  &  3 &  3 &  1 &  1 &  1 &  1 &  1 &  1   &&&&&\quad 2  &AR &I &N \\
\nex&\quad 6  &  3 &  2 &  2 &  1 &  1 &  1 &  1 &  1   &&&&&\quad 2  &AR &I &N \\
\nex&\quad 6  &  5 &  2 &  1 &  1 &  1 &  1 &  1   &&&&&&\quad 2  &AR &I &N \\
\nex&\quad 6  &  4 &  3 &  1 &  1 &  1 &  1 &  1   &&&&&&\quad 2  &AR &I &N \\
\nex&\quad 6  &  2 &  2 &  2 &  2 &  1 &  1 &  1 &  1   &&&&&\quad 2  &AR &I &N \\
\nex&\quad 6  &  4 &  2 &  2 &  1 &  1 &  1 &  1   &&&&&&\quad 2  &AR &I &N \\
\nex&\quad 6  &  3 &  3 &  2 &  1 &  1 &  1 &  1   &&&&&&\quad 2  &AR &I &N \\
\nex&\quad 6  &  5 &  3 &  1 &  1 &  1 &  1   &&&&&&&\quad 2  &AR &I &N \\
\nex&\quad 6  &  4 &  4 &  1 &  1 &  1 &  1   &&&&&&&\quad 2  &AR &I &N \\
\nex&\quad 6  &  3 &  2 &  2 &  2 &  1 &  1 &  1   &&&&&&\quad 2  &AR &I &N \\
\nex&\quad 6  &  5 &  2 &  2 &  1 &  1 &  1   &&&&&&&\quad 2  &AR &I &N \\
\nex&\quad 6  &  4 &  3 &  2 &  1 &  1 &  1   &&&&&&&\quad 2  &AR &I &N \\
\nex&\quad 6  &  3 &  3 &  3 &  1 &  1 &  1   &&&&&&&\quad 2  &AR &I &N \\
\nex&\quad 6  &  5 &  4 &  1 &  1 &  1   &&&&&&&&\quad 2  &AR &I &N \\
\nex&\quad 6  &  2 &  2 &  2 &  2 &  2 &  1 &  1   &&&&&&\quad 2  &AR &I &N \\
\nex&\quad 6  &  4 &  2 &  2 &  2 &  1 &  1   &&&&&&&\quad 2  &AR &I &N \\
\nex&\quad 6  &  3 &  3 &  2 &  2 &  1 &  1   &&&&&&&\quad 2  &AR &I &N \\
\nex&\quad 6  &  5 &  3 &  2 &  1 &  1   &&&&&&&&\quad 2  &AR &I &N \\
\nex&\quad 6  &  4 &  4 &  2 &  1 &  1   &&&&&&&&\quad 2  &AR &I &N \\
\nex&\quad 6  &  4 &  3 &  3 &  1 &  1   &&&&&&&&\quad 2  &AR &I &N \\
\nex&\quad 6  &  3 &  2 &  2 &  2 &  2 &  1   &&&&&&&\quad 2  &AR &P &N \\
\nex&\quad 6  &  5 &  2 &  2 &  2 &  1   &&&&&&&&\quad 2  &AR &P &N \\
\nex&\quad 6  &  4 &  3 &  2 &  2 &  1   &&&&&&&&\quad 2  &AR &P &N \\
\nex&\quad 6  &  3 &  3 &  3 &  2 &  1   &&&&&&&&\quad 2  &AR &P &N \\
\nex&\quad 6  &  3 &  3 &  2 &  2 &  2   &&&&&&&&\quad 2  &AR &P &N \\
\nex&\quad 8  &  3 &  3 &  3 &  3 &  3 &  1   &&&&&&&\quad 4  &AR &P &C \\
\nex&\quad 8  &  6 &  3 &  3 &  3 &  1   &&&&&&&&\quad 4  &AR &P &C \\
\nex&\quad 8  &  5 &  5 &  2 &  2 &  2   &&&&&&&&\quad 4  &AR &P &C \\
\nex&\quad 8  &  4 &  3 &  3 &  3 &  3   &&&&&&&&\quad 4  &AR &P &C \\
\tablerule
}}
\newpage

\vbox{\offinterlineskip\table
\def\marker{\vrule height 0.122truein depth 0.03truein width 0.7pt}
\def\topmarker{\vrule height 0.12truein depth 0.04truein width 0.7pt}
\def\tablerule{\noalign{\hrule height0.7pt}}
\let\\\cr
\def\s{\hglue 3.9pt}
\tabskip=0pt
\halign{\marker\hfil\s#\s\hfil&\marker\hfil\s#\s\hfil
&\marker\hfil\s#\s\hfil&\hfil\s#\s\hfil&\hfil\s#\s\hfil&\hfil\s#\s\hfil
&\hfil\s#\s\hfil&\hfil\s#\s\hfil&\hfil\s#\s\hfil&\hfil\s#\s\hfil
&\hfil\s#\s\hfil&\hfil\s#\s\hfil&\hfil\s#\s\hfil&\hfil\s#\s\hfil
&\marker\hfil\s#\s\hfil&\marker\hfil\s#\s\hfil
&\marker\hfil\s#\s\hfil&\marker\hfil\s#\s\hfil\marker\cr
\tablerule
\omit\topmarker&\omit\topmarker\s{\tiny Denominator}\s&\omit\topmarker\hfil
&\multispan 5
\hfil{\tiny Numerators}\hfil&&&&&&&\omit\topmarker\s{\tiny degree}\s
&\omit\topmarker\s{\tiny arithmetic?}\s
&\omit\topmarker\s{\tiny pure?}\s
&\omit\topmarker\s{\tiny compact?}\s\topmarker\cr
\tablerule
\nex&\quad 9  &  4 &  4 &  4 &  4 &  2   &&&&&&&&\quad 6  &AR &P &C \\
\nex&\quad 10  &  7 &  4 &  4 &  4 &  1   &&&&&&&&\quad 4  &AR &P &C \\
\nex&\quad 10  &  3 &  3 &  3 &  3 &  3 &  3 &  2   &&&&&&\quad 4  &AR &I &C \\
\nex&\quad 10  &  6 &  3 &  3 &  3 &  3 &  2   &&&&&&&\quad 4  &AR &I &C \\
\nex&\quad 10  &  9 &  3 &  3 &  3 &  2   &&&&&&&&\quad 4  &AR &I &C \\
\nex&\quad 10  &  6 &  6 &  3 &  3 &  2   &&&&&&&&\quad 4  &AR &I &C \\
\nex&\quad 10  &  5 &  3 &  3 &  3 &  3 &  3   &&&&&&&\quad 4  &AR &I &C \\
\nex&\quad 10  &  8 &  3 &  3 &  3 &  3   &&&&&&&&\quad 4  &AR &I &C \\
\nex&\quad 10  &  6 &  5 &  3 &  3 &  3   &&&&&&&&\quad 4  &AR &I &C \\
\nex&\quad 12  &  8 &  5 &  5 &  5 &  1   &&&&&&&&\quad 4  &AR &P &C \\
\nex&\quad 12  &  7 &  7 &  2 &  2 &  2 &  2 &  2   &&&&&&\quad 4  &AR &I &C \\
\nex&\quad 12  &  9 &  7 &  2 &  2 &  2 &  2   &&&&&&&\quad 4  &AR &I &C \\
\nex&\quad 12  &  7 &  7 &  4 &  2 &  2 &  2   &&&&&&&\quad 4  &AR &I &C \\
\nex&\quad 12 &  11 &  7 &  2 &  2 &  2   &&&&&&&&\quad 4  &AR &I &C \\
\nex&\quad 12  &  9 &  9 &  2 &  2 &  2   &&&&&&&&\quad 4  &AR &I &C \\
\nex&\quad 12  &  9 &  7 &  4 &  2 &  2   &&&&&&&&\quad 4  &AR &I &C \\
\nex&\quad 12  &  7 &  7 &  6 &  2 &  2   &&&&&&&&\quad 4  &AR &I &C \\
\nex&\quad 12  &  7 &  7 &  4 &  4 &  2   &&&&&&&&\quad 4  &AR &P &C \\
\nex&\quad 12  &  7 &  5 &  3 &  3 &  3 &  3   &&&&&&&\quad 4  &NA &P &N \\
\nex&\quad 12  &  5 &  5 &  5 &  3 &  3 &  3   &&&&&&&\quad 4  &AR &P &C \\
\nex&\quad 12 &  10 &  5 &  3 &  3 &  3   &&&&&&&&\quad 4  &AR &P &C \\
\nex&\quad 12  &  8 &  7 &  3 &  3 &  3   &&&&&&&&\quad 4  &NA &P &C \\
\nex&\quad 12  &  8 &  5 &  5 &  3 &  3   &&&&&&&&\quad 4  &AR &P &C \\
\nex&\quad 12  &  7 &  6 &  5 &  3 &  3   &&&&&&&&\quad 4  &NA &P &N \\
\nex&\quad 12  &  6 &  5 &  5 &  5 &  3   &&&&&&&&\quad 4  &AR &P &C \\
\nex&\quad 12  &  7 &  5 &  4 &  4 &  4   &&&&&&&&\quad 4  &NA &P &N \\
\nex&\quad 12  &  6 &  5 &  5 &  4 &  4   &&&&&&&&\quad 4  &NA &P &C \\
\nex&\quad 12  &  5 &  5 &  5 &  5 &  4   &&&&&&&&\quad 4  &AR &P &C \\
\nex&\quad 14 &  11 &  5 &  5 &  5 &  2   &&&&&&&&\quad 6  &AR &I &C \\
\nex&\quad 14  &  8 &  5 &  5 &  5 &  5   &&&&&&&&\quad 6  &AR &I &C \\
\nex&\quad 15  &  8 &  6 &  6 &  6 &  4   &&&&&&&&\quad 8  &NA &P &C \\
\nex&\quad 18 &  11 &  8 &  8 &  8 &  1   &&&&&&&&\quad 6  &AR &P &C \\
\nex&\quad 18 &  13 &  7 &  7 &  7 &  2   &&&&&&&&\quad 6  &NA &I &C \\
\nex&\quad 18 &  10 & 10 &  7 &  7 &  2   &&&&&&&&\quad 6  &AR &I &C \\
\nex&\quad 18 &  14 & 13 &  3 &  3 &  3   &&&&&&&&\quad 6  &AR &I &C \\
\nex&\quad 18 &  10 &  7 &  7 &  7 &  5   &&&&&&&&\quad 6  &AR &I &C \\
\nex&\quad 18  &  8 &  7 &  7 &  7 &  7   &&&&&&&&\quad 6  &NA &I &C \\
\nex&\quad 20 &  14 & 11 &  5 &  5 &  5   &&&&&&&&\quad 8  &NA &P &C \\
\nex&\quad 20 &  13 &  9 &  6 &  6 &  6   &&&&&&&&\quad 8  &NA &I &C \\
\nex&\quad 20 &  10 &  9 &  9 &  6 &  6   &&&&&&&&\quad 8  &NA &I &C \\
\nex&\quad 24 &  19 & 17 &  4 &  4 &  4   &&&&&&&&\quad 8  &NA &I &C \\
\nex&\quad 24 &  14 &  9 &  9 &  9 &  7   &&&&&&&&\quad 8  &NA &P &C \\
\nex&\quad 30 &  26 & 19 &  5 &  5 &  5   &&&&&&&&\quad 8  &AR &I &C \\
\nex&\quad 30 &  23 & 22 &  5 &  5 &  5   &&&&&&&&\quad 8  &NA &I &C \\
\nex&\quad 30 &  22 & 11 &  9 &  9 &  9   &&&&&&&&\quad 8  &AR &I &C \\
\nex&\quad 42 &  34 & 29 &  7 &  7 &  7   &&&&&&&&\quad 12  &NA &I &C \\
\nex&\quad 42 &  26 & 15 & 15 & 15 & 13   &&&&&&&&\quad 12  &NA &I &C \\
\tablerule}}


\begin{thebibliography}

\bibitem{Cooper-Thurston}
{\bf D~Cooper}, {\bf W\,P Thurston},
{\em Triangulating $3$--manifolds using $5$ vertex links},
Topology, 27 (1988) 23--25

\bibitem{Deligne-Mostow:Monodromy}
{\bf P~Deligne}, {\bf G\,D Mostow},
{\em Monodromy of hypergeometric functions and nonlattice integral
  monodromy},
Inst. Hautes \'Etudes Sci. Publ. Math. 63 (1986) 5--89

\bibitem{Deligne-Mostow:Commensurabilities}
{\bf P~Deligne}, {\bf G\,D Mostow},
{\em Commensurabilities among lattices in ${\rm {P}{U}}(1,n)$},
Annals of Mathematics Studies,  132,
Princeton University Press, Princeton, NJ (1993)

\bibitem{MR89c:52014}
{\bf D\,B\,A Epstein}, {\bf A~Marden},
{\em Convex hulls in hyperbolic space, a theorem of {S}ullivan, and
  measured pleated surfaces},
from: ``Analytical and geometric aspects of hyperbolic space
  (Coventry/Durham, 1984)'', Cambridge Univ. Press, Cambridge
  (1987) 113--253

\bibitem{Deligne-Mostow:Picard}
{\bf G\,D Mostow},
{\em Generalized {P}icard lattices arising from half-integral conditions}
Inst. Hautes \'Etudes Sci. Publ. Math. 63 (1986) 91--106

\bibitem{Picard:hypergeometriques-2}
{\bf E~Picard},
{\em Sur les fonctions hyperfuchsiennes provenant des s\'eries
  hyperg\'eometriques de deux variables},
Ann. ENS III, 2 (1885) 357--384

\bibitem{Picard:hypergeometriques-1}
{\bf E~Picard},
{\em Sur une extension aux fonctions de deux variables du probl\`eme de
  Riemann relatif aux fonctions hyperg\'eometriques},
Bull. Soc. Math. Fr. 15 (1887) 148--152

\bibitem{MR96h:57010}
{\bf Igor Rivin},
{\em Euclidean structures on simplicial surfaces and hyperbolic volume},
Annals of Math. 139 (1994) 553--580

\bibitem{MR96i:52008}
{\bf Igor Rivin},
{\em A characterization of ideal polyhedra in hyperbolic $3$--space},
Annals of Math. 143 (1996) 51--70

\bibitem{MR93j:52015}
{\bf Igor Rivin}, {\bf Craig~D Hodgson}, 
{\em A characterization of compact convex polyhedra in hyperbolic
  $3$--space},
Invent. Math. 111 (1993) 77--111

\bibitem{Thurston:GT3M}
{\bf William~P Thurston},
{\em Geometry and Topology of Three--Manifolds},
Princeton lecture notes (1979)\qua
{\tt  http://www.msri.org/publications/books/gt3m}

\bibitem{MR91h:53059}
{\bf Marc Troyanov},
{\em Prescribing curvature on compact surfaces with conical singularities},
Trans. Amer. Math. Soc. 324 (1991) 793--821

\end{thebibliography}
\end{document}